%% file: reglie.tex
\input amstex
\input amsppt.sty   
\input diag.tex

\magnification=\magstep1
\hsize 30pc
\vsize 47pc
\def\nmb#1#2{#2}         
\def\cit#1#2{\ifx#1!\cite{#2}\else#2\fi} 
\def\totoc{}             
\def\idx{}               
\def\ign#1{}             

\redefine\o{\circ}
\define\X{\frak X}
\define\al{\alpha}
\define\be{\beta}
\define\ga{\gamma}
\define\de{\delta}
\define\ep{\varepsilon}
\define\ze{\zeta}
\define\et{\eta}
\define\th{\theta}

\define\ka{\kappa}
\define\la{\lambda}

\define\si{\sigma}
\define\ta{\tau}
\define\ph{\varphi}

\define\ps{\psi}
\define\om{\omega}

\define\Th{\Theta}

\define\Ph{\Phi}
\define\Ps{\Psi}
\define\Om{\Omega}
\define\x{\times}

\define\Fl{\operatorname{Fl}}

\define\sign{\operatorname{sign}}

\define\conj{\operatorname{conj}}

\define\ev{\operatorname{ev}}
\define\Diff{\operatorname{Diff}}

\define\Pt{\operatorname{Pt}}

\define\Id{\operatorname{Id}}
\redefine\i{^{-1}}
\define\row#1#2#3{#1_{#2},\ldots,#1_{#3}}
\define\rowup#1#2#3{#1^{#2},\ldots,#1^{#3}}
\predefine\LL{\L}
\redefine\L{{\Cal L}}

\define\ddt{\left.\tfrac \partial{\partial t}\right\vert_0}
\define\dd#1{\tfrac \partial{\partial #1}}
\define\pr{\operatorname{pr}}


\define\ins{{\operatorname{ins}}}
\define\ad{{\operatorname{ad}}}
\define\Ad{{\operatorname{Ad}}}

\define\g{{\frak g}}
\define\h{{\frak h}}

\define\evol{{\operatorname{evol}}}
\define\Evol{{\operatorname{Evol}}}
\def\today{\ifcase\month\or
 January\or February\or March\or April\or May\or June\or
 July\or August\or September\or October\or November\or December\fi
 \space\number\day, \number\year}
\topmatter
\title  Regular infinite dimensional Lie groups
\endtitle
\author Andreas Kriegl\\ Peter W. Michor  \endauthor
\affil
Institut f\"ur Mathematik, Universit\"at Wien,\\
Strudlhofgasse 4, A-1090 Wien, Austria.\\
Erwin Schr\"odinger Institut f\"ur Mathematische Physik,
Pasteurgasse 6/7, A-1090 Wien, Austria
\endaffil
\address
Institut f\"ur Mathematik, Universit\"at Wien,
Strudlhofgasse 4, A-1090 Wien, \newline Austria
\endaddress
\email kriegl\@pap.univie.ac.at 
\endemail
\address
Erwin Schr\"odinger Institut f\"ur Mathematische Physik,
Pasteurgasse 6/7, A-1090 Wien, Austria
\endaddress
\email Peter.Michor\@esi.ac.at \endemail
\date {\today} \enddate
\abstract 
Regular Lie groups are infinite dimensional Lie groups with the 
property that smooth curves in the Lie algebra integrate to 
smooth curves in the group in a smooth way (an `evolution operator' 
exists). Up to now all known smooth Lie groups are regular.
We show in this paper that 
regular Lie groups allow to push surprisingly far the geometry of 
principal bundles: parallel transport exists and flat connections 
integrate to horizontal foliations as in finite dimensions. 
As consequences we obtain that Lie algebra homomorphisms intergrate
to Lie group homomorphisms, if the source group is simply connected 
and the image group is regular.
\endabstract
\endtopmatter

\document

\heading Table of contents \endheading
\noindent 1. Introduction 
\leaders \hbox to 1em{\hss .\hss }\hfill {\eightrm 1}\par 
\noindent 2. Calculus of smooth mappings 
\leaders \hbox to 1em{\hss .\hss }\hfill {\eightrm 2}\par 
\noindent 3. Lie groups 
\leaders \hbox to 1em{\hss .\hss }\hfill {\eightrm 7}\par 
\noindent 4. Bundles and connections
\leaders \hbox to 1em{\hss .\hss }\hfill {\eightrm 12}\par 
\noindent 5. Regular Lie groups 
\leaders \hbox to 1em{\hss .\hss }\hfill {\eightrm 16}\par 
\noindent 6. Bundles with regular structure groups
\leaders \hbox to 1em{\hss .\hss }\hfill {\eightrm 28}\par 
\noindent 7. Rudiments of Lie theory for regular Lie groups 
\leaders \hbox to 1em{\hss .\hss }\hfill {\eightrm 32}\par 

\head\totoc\nmb0{1}. Introduction \endhead

On the one hand the theory of infinite dimensional Lie groups and Lie 
algebras is very rich: Kac-Moody algebras and the Virasoro algebra 
have a rich and important theory of representations and many 
applications, and subgroups of diffeomorphism groups play an 
extremely important role in differential topology, differential 
geometry, and general relativity. On the other hand classical Lie 
theory carries over to them only in rare pieces: There are (even 
Banach) Lie algebras without Lie groups, see \cit!{3} and \cit!{7}, 
and the exponential mapping in general is not surjective onto any 
neighborhood of the identity. The most surprising result in this 
direction is \cit!{5}, where it is shown that in the diffeomorphism 
group of any manifold of dimension at least 2 one can find a smooth 
curve of diffeomorphisms starting at the identity such that the 
points of this curve form a set of generators for a free subgroup of 
the diffeomorphism group which meets the image of the exponential 
mapping only in the identity.

In view of these difficulties 
the theory of infinite dimensional Lie groups and Lie 
can be pushed surprisingly far: Exponential mappings are unique if they 
exist, and then one can give a formula for their derivatives, see 
\cit!{6} and \nmb!{5.9} below. 

In \cit!{14} and \cit!{15} the notion of `regular Fr\'echet Lie 
group' was introduced in an attempt to find conditions which ensure 
the existence of exponential mappings: certain product integrals 
were required to converge. Their main result was that the invertible 
Fourier integral operators form a regular Fr\'echet Lie group with 
the space of pseudo differential operators as Lie algebra, see 
\cit!{17}, and also \cit!{1} for a more general group of Fourier integral 
operators, but without regularity.
In \cit!{13} Milnor weakened this to the 
assumption that smooth curves in the Lie algebra integrate to 
smooth curves in the group in a smooth way (an `evolution operator' 
exists), and it is this notion which we take up in this paper, but for 
general Lie groups modelled on locally convex spaces, where we use 
the convenient calculus from \cit!{4}. 
Up to now nobody found a non-regular Lie group.

We show in this paper that the use of 
regular Lie groups allows to push surprisingly far the geometry of 
principal bundles: parallel transport exists and flat connections 
integrate to horizontal foliations as in finite dimensions. 
As consequences we obtain that Lie algebra homomorphisms intergrate
to Lie group homomorphisms, if the source group is simply connected 
and the image group is regular.

The actual development is quite involved. We start with general 
infinite dimensional Lie groups in section \nmb!{3}, 
but for a detailed study of the 
evolution operator of regular Lie groups (\nmb!{5.3}) we need in 
\nmb!{5.9} the Maurer-Cartan 
equation for right (or left) logarithmic derivatives (\nmb!{5.1}) 
of mappings with values in the Lie group, 
and this we can get only by looking at 
principal connections. Thus section \nmb!{4} treats principal 
bundles, connections, and curvature as far as we shall need them.
We can then prove the strong existence results mentioned 
above and treat regular Lie groups in section \nmb!{5}, and 
principal bundles with regular structure groups in section 
\nmb!{6}. The last section \nmb!{7} develops rudiments of Lie 
theory for regular Lie groups as sketched above.

These results were obtained in a systematic study of properties of 
regular Lie groups for the book in preparation \cit!{11}, where also 
many of the known Lie groups are treated and are shown to be regular.

\heading\totoc\nmb0{2}. Calculus of smooth mappings \endheading

The traditional differential calculus works 
well for finite dimensional vector spaces and for Banach spaces. For 
more general locally convex spaces a whole flock of different 
theories were developed, most of them rather complicated and not 
really convincing. The main difficulty is that the composition of 
linear mappings stops to be jointly continuous at the level of Banach 
spaces, for any compatible topology. 

We shall use in this paper the calculus in infinite dimensions as 
developed in \cit!{4}.

\subheading{\nmb.{2.1}. Convenient vector spaces} Let $E$ be a 
locally convex vector space. 
A curve $c:\Bbb R\to E$ is called 
{\it smooth} or $C^\infty$ if all derivatives exist and are 
continuous - this is a concept without problems. Let 
$C^\infty(\Bbb R,E)$ be the space of smooth functions. It can be 
shown that $C^\infty(\Bbb R,E)$ does not depend on the locally convex 
topology of $E$, but only on its associated bornology (system of bounded 
sets).

$E$ is said to be a {\it convenient 
vector space} if one of the following equivalent
conditions is satisfied (called $c^\infty$-completeness):
\roster
\item For any $c\in C^\infty(\Bbb R,E)$ the (Riemann-) integral 
       $\int_0^1c(t)dt$ exists in $E$.
\item A curve $c:\Bbb R\to E$ is smooth if and only if $\la\o c$ is 
       smooth for all $\la\in E'$, where $E'$ is the dual consisting 
       of all continuous linear functionals on $E$.
\item Any Mackey-Cauchy-sequence (i\.e\. $t_{nm}(x_n-x_m)\to 0$  
       for some $t_{nm}\to \infty$ in $\Bbb R$) converges in $E$. 
       This is visibly a weak completeness requirement.
\endroster

The final topology with respect to all smooth curves is called the 
$c^\infty$-topology on $E$, which then is denoted by $c^\infty E$. 
For Fr\'echet spaces it coincides with 
the given locally convex topology, but on the space $\Cal D$ of test 
functions with compact support on $\Bbb R$ it is strictly finer.

\subheading{\nmb.{2.2}. Smooth mappings} Let $E$ and $F$ be locally 
convex vector spaces, and let $U\subset E$ be $c^\infty$-open. 
A mapping $f:U\to F$ is called {\it smooth} or 
$C^\infty$, if $f\o c\in C^\infty(\Bbb R,F)$ for all 
$c\in C^\infty(\Bbb R,U)$.

\proclaim{\nmb.{2.3}. Results}
The main properties of smooth calculus are the following.
\roster
\item For mappings on Fr\'echet spaces this notion of smoothness 
       coincides with all other reasonable definitions. Even on 
       $\Bbb R^2$ this is non-trivial, see \cit!{2}.
\item Multilinear mappings are smooth if and only if they are 
       bounded.
\item If $f:E\supseteq U\to F$ is smooth then the derivative 
       $df:U\x E\to F$ is  
       smooth, and also $df:U\to L(E,F)$ is smooth where $L(E,F)$ 
       denotes the space of all bounded linear mappings with the 
       topology of uniform convergence on bounded subsets.
\item The chain rule holds.
\item The space $C^\infty(U,F)$ is again a convenient vector space 
       where the structure is given by the obvious injection
$$
C^\infty(U,F)\to \prod_{c\in C^\infty(\Bbb R,U)} C^\infty(\Bbb R,F)
\to \prod \Sb c\in C^\infty(\Bbb R,U) \\ \la\in F' 
\endSb C^\infty(\Bbb R,\Bbb R).
$$
\item The exponential law holds:
$$
C^\infty(U,C^\infty(V,G)) \cong C^\infty(U\x V, G)
$$
     is a linear diffeomeorphism of convenient vector spaces. Note 
     that this is the main assumption of variational calculus.
\item A linear mapping $f:E\to C^\infty(V,G)$ is smooth (bounded) if 
       and only if $E @>f>> C^\infty(V,G) @>{\ev_v}>> G$ is smooth 
       for each $v\in V$. This is called the smooth uniform 
       boundedness theorem and it is quite applicable.
\endroster
\endproclaim

\example{\nmb.{2.4} Counterexamples in infinite dimensions against 
common beliefs on ordinary differential equations}
Let $E:=s$ be the Fr\'echet space of rapidely decreasing sequence 
(Note that by the theory of Fourier series we have $s=C^\infty(S^1,\Bbb R)$) 
and
consider the continuous linear operator $T:E\to E$ given by
$T(x_0,x_1,x_2,\dots):=(0,1^2x_1,2^2x_2,3^2x_3,\dots)$. 
The ordinary linear differential equation $x'(t)=T(x(t))$ 
with constant coefficients has no solution in $s$ for certain initial 
values. By recursion one sees that 
the general solution should be given by 
$$
x_n(t)=\sum_{i=0}^n
     \left(\tfrac{n!}{i!}\right)^2 x_i(0)\frac{t^{n-i}}{(n-i)!}
$$
If the initial value is a finite sequence, say $x_n(0)=0$ for $n>N$
and $x_N(0)\ne 0$, then 
$$\align
x_n(t)&=\sum_{i=0}^N
     \left(\tfrac{n!}{i!}\right)^2 x_i(0)\frac{t^{n-i}}{(n-i)!}\\
&=\frac{(n!)^2}{(n-N)!}\;t^{n-N}\sum_{i=0}^N
     \left(\tfrac{1}{i!}\right)^2 x_i(0)\tfrac{(n-N)!}{(n-i)!}\;t^{N-i}\\
|x_n(t)| &\ge \frac{(n!)^2}{(n-N)!}\;|t|^{n-N}
     \left(|x_N(0)|\left(\tfrac1{N!}\right)^2 -
          \sum_{i=0}^{N-1}
     \left(\tfrac{1}{i!}\right)^2 
     |x_i(0)|\tfrac{(n-N)!}{(n-i)!}|t|^{N-i}
     \right)\\
&\ge \frac{(n!)^2}{(n-N)!}\;|t|^{n-N}
     \left(|x_N(0)|\left(\tfrac1{N!}\right)^2 -
          \sum_{i=0}^{N-1}
     \left(\tfrac{1}{i!}\right)^2 |x_i(0)||t|^{N-i}
     \right)\\
\endalign$$
where the first factor does not lie in the space $s$ of rapidly 
decreasing sequences and where the second factor is larger than 
$\ep>0$ for $t$ small enough.
So at least for a dense set of initial values this differential 
equation has no local solution. 

This shows also, that the theorem of 
Frobenius is wrong, in the 
following sense: The vector field $x\mapsto T(x)$ generates a 
1-dimensional subbundle $E$ of the tangent bundle on the open subset 
$s\setminus 0$. It is involutive since it is 1-dimensional. But 
through points representing finite sequences there exist no local integral 
submanifolds ($M$ with $TM=E|M$). Namely, if $c$ were a smooth 
nonconstant curve with $c'(t)=f(t).T(c(t))$ for some smooth 
function $f$, then $x(t):=c(h(t))$ would satisfy $x'(t)=T(x(t))$, where 
$h$ is a solution of $h'(t)=1/f(h(t))$. 

As next example consider 
$E:=\Bbb R^\Bbb N$ and the continuous linear operator 
$T:E\to E$ given by
$T(x_0,x_1,\dots):= (x_1,x_2,\dots )$.
The corresponding differential equation has solutions for every 
initial value $x(0)$, since the coordinates must satisfy the recusive 
relations $x_{k+1}(t)=x_k'(t)$ and hence any smooth functions 
$x_0:\Bbb R\to \Bbb R$ gives rise to a
solution $x(t):=(x_0^{(k)}(t))_k$ with initial value
$x(0)=(x_0^{(k)}(0))_k$. So by Borel's theorem there exist solutions 
to this equation for any initial value and the difference of any two 
functions with same initial value
is an arbitray infinite flat function.
Thus the solutions are far from being unique.
Note that $\Bbb R^\Bbb N$ is a topological direct summand in 
$C^\infty(\Bbb R,\Bbb R)$
via the projection $f\mapsto (f(n))_n$, and hence the same situation
occurs in $C^\infty(\Bbb R,\Bbb R)$.

Let now $E:=C^\infty(\Bbb R,\Bbb R)$ and consider the continuous 
linear operator
$T:E\to E$ given by $T(x):=x'$. Let 
$x:\Bbb R\to C^\infty(\Bbb R,\Bbb R)$ be a
solution of the equation $x'(t)=T(x(t))$. In terms of 
$\hat x:\Bbb R^2\to \Bbb R$
this says 
$\frac{\partial}{\partial t}\hat x(t,s)
     =\frac{\partial}{\partial s}\hat x(t,s)$. 
Hence
$r\mapsto \hat x(t-r,s+r)$ has vanishing derivative everywhere and so 
this function is constant, and in particular 
$x(t)(s)=\hat x(t,s)=\hat x(0,s+t)=x(0)(s+t)$.
Thus we have a smooth solution $x$ uniquely determined by the initial
value $x(0)\in C^\infty(\Bbb R,\Bbb R)$ which even describes a flow 
for the vector field $T$ in the sense of \nmb!{2.7} below.
In general this solution is however not real-analytic, since for any
$x(0)\in C^\infty(\Bbb R,\Bbb R)$, which is not real-analytic in a 
neighborhood of
a point $s$ the composite $\ev_s\o x=x(s+\quad)$ is not real-analytic 
around $0$.
\endexample

\subhead\nmb.{2.5}. Manifolds \endsubhead
In the sequel we shall use smooth manifolds $M$ modelled on 
$c^\infty$-open subsets of convenient vector spaces. See 
\cit!{10} for an account of this. Since we shall need it we also 
include some results on vector fields and their flows.

\proclaim{\nmb.{2.6}. Lemma} 
Consider vector fields $X_i\in C^\infty(TM)$
and $Y_i\in C^\infty(TN)$ for $i=1,2$, and a smooth mapping $f:M\to N$.
If $X_i$ and $Y_i$ are $f$-related for $i=1,2$, i\.e\. $Tf\o X_i = 
Y_i\o f$, then also
$[X_1,X_2]$ and $[Y_1,Y_2]$ are $f$-related.
\endproclaim

\demo{Proof} 
We choose $h\in C^\infty(N,\Bbb R)$ and we view each vector field as 
a derivation. This is possible if we either have smooth partitions of 
unity or if we pass to sheaves of smooth functions. The converse is 
wrong in general, see \cit!{10} and \cit!{11}. 
Then by assumption we have 
$T f\o X_i = Y_i\o f$, thus:
$$\multline (X_i(h\o f))(x) = X_i(x)(h\o f) = (T_xf.X_i(x))(h) =\\
=(T f\o X_i)(x)(h) = (Y_i\o f)(x)(h) = Y_i(f(x))(h) = (Y_i(h))(f(x)),
\endmultline$$
so $X_i(h\o f) = (Y_i(h))\o f$, and we may continue:
$$\align [X_1,X_2](h\o f) &= X_1(X_2(h\o f)) - X_2(X_1(h\o f))=\\ 
&= X_1(Y_2(h)\o f) - X_2(Y_1(h)\o f) = \\
&= Y_1(Y_2(h))\o f - Y_2(Y_1(h))\o f = [Y_1,Y_2](h)\o f. 
\endalign$$
But this means  $T f\o [X_1,X_2] = [Y_1,Y_2]\o f$.
\qed\enddemo

In particular if $f:M\to N$ is a local
diffeomorphism (so $(T_xf)\i$ makes sense for each $x\in M$),
then for $Y\in C^\infty(T N)$ a vector field 
$f^*Y\in C^\infty(T M)$ is defined
by $(f^*Y)(x) = (T_xf)\i.Y(f(x))$. The linear mapping
$f^*:C^\infty(T N) \to C^\infty(T M)$ is then a Lie algebra 
homomorphism.

\subhead\nmb.{2.7}. The flow of a vector field \endsubhead
Let $X\in C^\infty(TM)$ be a vector field. A 
\idx{\it local flow} 
$\Fl^X$ for $X$ is a smooth mapping 
$\Fl^X: M\x\Bbb R\supset U \to M$ defined on a $c^\infty$-open 
neighborhood $U$ of $M\x0$ such that
\roster
\item $\tfrac d{dt} \Fl^X_t(x)=X(\Fl^X_t(x))$.
\item $\Fl^X_0(x)=x$ for all $x\in M$.
\item $U\cap (\{x\}\x \Bbb R)$ is a connected open interval. 
\item $\Fl^X_{t+s} = \Fl^X_t\o \Fl^X_s$ holds in the following sense. 
     If the right hand side exists then also the left hand side 
     exists and we have equality. Moreover:
     If $\Fl^X_s$ exists, then the existence of both 
     sides is equivalent and they are equal.
\endroster

\proclaim{\nmb.{2.8}. Lemma} Let $X\in C^\infty(TM)$ be a 
vector field which admits a local flow $\Fl^X_t$. 
Then each for integral curve $c$ of $X$ we have $c(t)=\Fl^X_t(c(0))$, 
thus there exists a unique maximal flow. 
Furthermore $X$ is $\Fl^X_t$-related to itself, i\.e\. 
$T(\Fl^X_t)\o X = X\o \Fl^X_t$. 
\endproclaim

\demo{Proof} We compute 
$$\align
\tfrac d{dt} \Fl^X(-t,c(t)) &= -\tfrac d{ds}|_{s=-t}\Fl^X(s,c(t)) + 
     \tfrac d{ds}|_{s=t} \Fl^X(-t,c(s))\\
&= -\tfrac d{ds}|_{s=0}\Fl^X_{-t}\Fl^X(s,c(t)) + T(\Fl^X_{-t}).c'(t)\\
&= -T(\Fl^X_{-t}).X(c(t)) + T(\Fl^X_{-t}).X(c(t)) =0,\\
\endalign$$
Thus $\Fl^X_{-t}(c(t))=c(0)$ is constant, so $c(t)=\Fl^X_t(c(0))$.
For the second assertion we have 
$X\o \Fl^X_t = \tfrac d{dt}\Fl^X_t = \tfrac d{ds}|_0 \Fl^X_{t+s} = 
\tfrac d{ds}|_0 (\Fl^X_t\o\Fl^X_s) = T(\Fl^X_t)\o\tfrac d{ds}|_0\Fl^X_s 
= T(\Fl^X_t)\o X$. 
\qed\enddemo

\proclaim{\nmb.{2.9}. Lemma} Let $X\in C^\infty(TM)$ and 
$Y\in C^\infty(TN)$ be $f$-related vector fields for a smooth mapping 
$f:M \to N$ which have local flows $\Fl^X$ and $\Fl^Y$. 
Then we have $f\o \Fl^X_t = \Fl^Y_t\o f$, whenever both
sides are defined. 

Moreover, if $f$ is a diffeomorphism we
have $\Fl^{f^*Y}_t = f\i\o \Fl^Y_t\o f$ in the following sense: If 
one side exists then also the other and they are equal. 
\endproclaim
For $f=Id_M$ this again implies that if there exists a flow then 
there exists a unique maximal flow $\Fl^X_t$.

\demo{Proof} We have $Y\o f = Tf\o X$ and thus (using \nmb!{2.7}.3 
and \nmb!{2.8}) for small $t$ we get
$$\align
\tfrac d{dt}(\Fl^Y_t\o f\o \Fl^X_{-t}) 
&= Y\o \Fl^Y_t\o f\o \Fl^X_{-t} - T(\Fl^Y_t)\o Tf\o X\o \Fl^X_{-t} \\
&= T(\Fl^Y_t)\o Y\o f\o \Fl^X_{-t} - T(\Fl^Y_t)\o Tf\o X\o \Fl^X_{-t} 
= 0.
\endalign$$ 
So $(\Fl^Y_t\o f\o \Fl^X_{-t})(x) = f(x)$ or 
$f(\Fl^X_t(x)) = \Fl^Y_t(f(x))$ for small $t$. 
By the flow properties (\nmb!{2.7}.4) we get the result by a 
connectedness argument as follows:
In the common interval of definition we consider the closed subset
$J_x:=\{t: f(\Fl^X_t(x)) = \Fl^Y_t(f(x))\}$. This set is also open since 
for $t\in J_x$ and small $|s|$ we have
$f(\Fl^X_{t+s}(x)) = f(\Fl^X_s(\Fl^X_t(x))) = \Fl^Y_s(f(\Fl^X_t(x))) 
= \Fl^Y_s(\Fl^Y_t(f(x))) = \Fl^Y_{t+s}(f(x))$.
\qed\enddemo

\subhead \nmb.{2.10}. The Lie derivative \endsubhead
We will meet situations (in \nmb!{4.2}) where we do not know that 
the flow of $X$ exists but where we will be able to produce the 
following assumption: 
Suppose that 
$\ph:\Bbb R\x M\supset U\to M$ is a smooth mapping such that 
$(t,x)\mapsto (t,\ph(t,x)=\ph_t(x))$ is a diffeomorphism 
$U\to V$, where $U$ and $V$ are open neighborhoods of 
$\{0\}\x M$ in $\Bbb R\x M$, and such that
$\ph_0=\Id_M$ and $\ddt \ph_t=X\in C^\infty(TM)$.
Then again 
$\tfrac d{dt}|_0(\ph_t)^*f = \tfrac d{dt}|_0f\o \ph_t = 
df\o X = X(f)$.

\proclaim{Lemma} In this situation we have for $Y\in C^\infty(TM)$, 
and for a $k$-form $\om\in \Om^k(M)$:
$$\gather
\tfrac d{dt}|_0(\ph_t)^*Y = [X,Y],\\
\dd t|_0 (\ph_t)^*\om = \L_X\om\\
\endgather$$
\endproclaim

\demo{Proof} Let $f\in C^\infty(M,\Bbb R)$ be a 
function and consider the mapping $\al(t,s):= Y(\ph(t,x))(f\o
\ph_s)$, which is 
locally defined near 0. It satisfies
$$\align 
&\al(t,0)=Y(\ph(t,x))(f), \\
&\al(0,s) = Y(x)(f\o\ph_s), \\
&\dd t\al(0,0) = \ddt Y(\ph(t,x))(f) = \ddt(Yf)(\ph(t,x)) =
        X(x)(Yf),\\ 
&\dd s\al(0,0) = \dd s|_0 Y(x)(f\o \ph_s) = Y(x)\dd s|_0(f\o
	\ph_s) = Y(x)(Xf).\endalign$$
So $\dd u|_0 \al(u,-u)=[X,Y]_x(f)$.
But on the other hand we have 
$$\align 
\dd u|_0 \al(u,-u) &= \dd u|_0 Y(\ph(u,x))(f\o\ph_{-u}) = \\
&= \dd u|_0\left(T(\ph_{-u})\o Y\o \ph_u\right)_x(f) \\
&= (\tfrac d{dt}|_0(\ph_t)^*Y)_x(f).
\endalign$$
We may identify $k$-forms on $M$ with 
$C^\infty(M,\Bbb R)$-multilinear mappings on vector fields (if smooth 
partitions of unity exist or if we pass to sheaves of vector fields). 
The converse is wrong in general, see \cit!{11}.
For (local) vector fields $Y_i\in C^\infty(TM)$ we have
$$\align
(\dd t|_0 (\ph_t)^*\om)&(Y_1,\dots,Y_k) 
= \dd t|_0 (\om((\ph_{-t})^*Y_1,\dots,(\ph_{-t})^*Y_k)\o \ph_{t})\\
&= \sum_{j=1}^k \om(Y_1,\dots,\dd t|_0 (\ph_{-t})^*Y_j,\dots,Y_k) 
+ \dd t|_0 (\ph_t)^*(\om(Y_1,\dots,Y_p) \\
&= X(\om(\row Y1k)) 
        - \sum_{i=1}^k \om(Y_1,\dots,[X,Y_i],\dots,Y_k)\\
&=\L_X\om(Y_1,\dots,Y_k).
\endalign$$
This is the usual formula for the Lie derivative.
\qed\enddemo

\head\totoc \nmb0{3}.  Lie groups \endhead

\subhead \nmb.{3.1}. Definition \endsubhead 
A \idx{\it Lie group} $G$ is a
smooth manifold modelled on $c^\infty$-open subsets of a convenient 
vector space, and a group such that the multiplication
$\mu:G\x G\to G$ and the inversion $\nu:G\to G$ are smooth. 
We shall use the following notation: \newline
$\mu:G\x G\to G$, multiplication, $\mu(x,y) = x.y$. \newline
$\mu_a: G\to G$, left translation, $\mu_a(x) = a.x$.\newline
$\mu^a: G\to G$, right translation, $\mu^a(x) = x.a$.\newline
$\nu: G\to G$, inversion, $\nu(x) = x\i$.\newline
$e\in G$, the unit element. 

\proclaim{\nmb.{3.2}. Lemma} The tangent mapping 
$T_{(a,b)}\mu:T_aG\x T_bG \to T_{ab}G$ is given by 
$$T_{(a,b)}\mu.(X_a,Y_b) = T_a(\mu^b).X_a + T_b(\mu_a).Y_b.$$
and $T_a\nu:T_aG\to T_{a\i}G$ is given by
$$T_a\nu = - T_e(\mu^{a\i}).T_a(\mu_{a\i}) 
= - T_e(\mu_{a\i}).T_a(\mu^{a\i}).$$ 
\endproclaim

\demo{Proof} Let $\ins_a:G\to G\x G$, $\ins_a(x)=(a,x)$ be the right
insertion and let $\ins^b:G\to G\x G$, $\ins^b(x)=(x,b)$ be the left
insertion. Then we have
$$\multline T_{(a,b)}\mu.(X_a,Y_b) 
= T_{(a,b)}\mu.(T_a(\ins^b).X_a + T_b(\ins_a).Y_b) = \\
= T_a(\mu\o \ins^b).X_a + T_b(\mu\o \ins_a).Y_b
= T_a(\mu^b).X_a + T_b(\mu_a).Y_b.
\endmultline$$
Now we differentiate the
equation $\mu(a,\nu(a))=e$; this gives in turn
$$\gather 0_e = T_{(a,a\i)}\mu.(X_a,T_a\nu.X_a) = T_a(\mu^{a\i}).X_a +
T_{a\i}(\mu_a).T_a\nu.X_a,\\
T_a\nu.X_a = - T_e(\mu_a)\i. T_a(\mu^{a\i}).X_a 
= - T_e(\mu_{a\i}).T_a(\mu^{a\i}).X_a.\qed\endgather$$ 
\enddemo

\subhead \nmb.{3.3}. Invariant vector fields and Lie algebras \endsubhead
Let $G$ be a (real) Lie group. A vector field $\xi$ on $G$ is
called \idx{\it left invariant}, 
if $\mu_a^*\xi = \xi$ for all $a\in
G$, where $\mu_a^*\xi = T(\mu_{a\i})\o\xi\o\mu_a$.
Since we have $\mu_a^*[\xi,\et] =
[\mu_a^*\xi,\mu_a^*\et]$, the space $\X_L(G)$ of all left invariant
vector fields on $G$ is closed under the Lie bracket, so it is a
sub Lie algebra of $\X(G)$. 
Any left invariant vector field
$\xi$ is uniquely determined by $\xi(e)\in T_eG$, since
$\xi(a)=T_e(\mu_a).\xi(e)$. Thus the Lie algebra $\X_L(G)$ of
left invariant vector fields is linearly isomorphic to $T_eG$,
and on $T_eG$ the Lie bracket on $\X_L(G)$ induces a Lie algebra
structure, whose bracket is again denoted by $[\quad,\quad]$.
This Lie algebra will be denoted as usual by $\frak g$,
sometimes by $\operatorname{Lie}(G)$.

We will also give a name to the isomorphism with the space of left 
invariant vector fields: $L:\frak g\to \X_L(G)$, $X\mapsto L_X$, 
where $L_X(a)= T_e\mu_a.X$. Thus $[X,Y] = [L_X,L_Y](e)$.

Similarly a vector field $\eta$ on $G$ is called 
\idx{\it right invariant}, if
$(\mu^a)^*\eta =\eta$ for all $a\in G$. If $\xi$ is left
invariant, then $\nu^*\xi$ is right invariant.
The right invariant vector fields form a sub Lie algebra
$\X_R(G)$ of $\X(G)$, which is again linearly isomorphic to
$T_eG$ and induces the negative of the Lie algebra structure on $T_eG$. 
We will denote by
$R:\frak g = T_eG\to \X_R(G)$ the isomorphism discussed, which
is given by $R_X(a)=T_e(\mu^a).X$.

\proclaim{\nmb.{3.4}. Lemma} If $L_X$ is a left invariant vector
field and $R_Y$ is a right invariant one, then $[L_X,R_Y]=0$.
So if the flows of $L_X$ and $R_Y$ exist, they commute.
\endproclaim

\demo{Proof} We consider the vector field $0\x L_X\in \X(G\x G)$, given by 
$(0\x L_X)(a,b)= (0_a,L_X(b))$. Then $T_{(a,b)}\mu.(0_a,L_X(b))=
T_a\mu^b.0_a + T_b\mu_a.L_X(b) = L_X(ab)$, so $0\x L_X$ is
$\mu$-related to $L_X$.
Likewise $R_Y\x0$ is $\mu$-related to $R_Y$. But then 
$0=[0\x L_X,R_Y\x0]$ is $\mu$-related to $[L_X,R_Y]$ by \nmb!{2.6}. 
Since $\mu$ is surjective, $[L_X,R_Y]=0$ follows.
\qed\enddemo

\proclaim{\nmb.{3.5}. Lemma}
Let $\ph:G\to H$ be a smooth homomorphism of Lie
groups.
Then $\ph' := T_e\ph:\frak g=T_eG \to \frak
h=T_eH$ is a Lie algebra homomorphism.
\endproclaim

\demo{Proof} For $X\in \frak g$ and $x\in G$ we have 
$$\align 
T_x\ph.L_X(x) &= T_x\ph.T_e\mu_x.X = T_e(\ph\o\mu_x).X \\
&=T_e(\mu_{\ph(x)}\o\ph).X = T_e(\mu_{\ph(x)}).T_e\ph.X 
     = L_{\ph'(X)}(\ph(x)).
\endalign$$
So $L_X$ is $\ph$-related to $L_{\ph'(X)}$. By \nmb!{2.6} the field
$[L_X,L_Y] = L_{[X,Y]}$ is $\ph$-related to
$[L_{\ph'(X)},L_{\ph'(Y)}] = L_{[\ph'(X),\ph'(Y)]}$. So we have 
$T\ph\o L_{[X,Y]} = L_{[\ph'(X),\ph'(Y)]}\o\ph$. If we evaluate
this at $e$ the result follows.
\qed\enddemo

\subhead \nmb.{3.6}.  One parameter subgroups \endsubhead 
Let $G$ be a Lie
group with Lie algebra $\frak g$. A \idx{\it one parameter subgroup}
of $G$ is a Lie group homomorphism $\al:(\Bbb R,+) \to G$, i.e.
a smooth curve $\al$ in $G$ with $\al(s+t)=\al(s).\al(t)$, and hence 
$\al(0)=e$.  

Note that a smooth mapping 
$\be:(-\ep,\ep)\to G$ satisfying $\be(t)\be(s)=\be(t+s)$ for $|t|$, 
$|s|$, $|t+s|<\ep$ is the restriction of a one parameter subgroup.
Namely, choose $0<t_0<\ep/2$. Any $t\in\Bbb R$ can be uniquely 
written as $t=N.t_0+t'$ for $0\le t'<t_0$ and $N\in\Bbb Z$. Put
$\al(t)=\be(t_0)^N\be(t')$. The required properties are easy to 
check.

\proclaim{ Lemma} Let $\al:\Bbb R\to G$ be a smooth curve with
$\al(0)=e$. Let $X\in \frak g$. 
Then the following assertions are equivalent.
\roster
\item $\al$ is a one parameter subgroup with $X=\ddt \al(t)$.
\item $\al(t)$ is an integral curve of the left invariant vector 
       field $L_X$, and also an integral curve of the right invariant 
       vector field $R_X$. 
\item $\Fl^{L_X}(t,x):= x.\al(t)$ (or
     $\Fl^{L_X}_t=\mu^{\al(t)}$) is the (unique by \nmb!{2.9}) 
       global flow of $L_X$ in the sense of \nmb!{2.7}.
\item $\Fl^{R_X}(t,x):= \al(t).x$ (or
     $\Fl^{R_X}_t=\mu_{\al(t)}$) is the (unique)
       global flow of $R_X$.
\endroster
Moreover each of these properties determines $\al$ uniquely.
\endproclaim

\demo{Proof} \therosteritem1 $\Longrightarrow$ \therosteritem3.
We have $\tfrac d{dt}x.\al(t) = \tfrac d{ds}|_0x.\al(t+s) = \tfrac
d{ds}|_0 x.\al(t).\al(s) =
 \tfrac d{ds}|_0 \mu_{x.\al(t)}\al(s)
= T_e(\mu_{x.\al(t)}).\tfrac d{ds}|_0 \al(s) = L_X(x.\al(t))$.
Since it is obviously a flow, we have \therosteritem3.

\therosteritem3 $\Longleftrightarrow$ \therosteritem4.
We have 
$\Fl^{\nu^*\xi}_t = \nu\i\o \Fl^\xi_t\o\nu$ by \nmb!{2.9}. Therefore 
we have by \nmb!{3.3}
$$\align
(\Fl^{R_X}_t(x\i))\i &= (\nu\o\Fl^{R_X}_t\o\nu)(x)
	= \Fl^{\nu^*R_X}_t(x) \\
&= \Fl^{-L_X}_{t}(x) = \Fl^{L_X}_{-t}(x) = x.\al(-t).
\endalign$$
So $\Fl^{R_X}_t(x\i) = \al(t).x\i$, and $\Fl^{R_X}_t(y) = \al(t).y$.

\therosteritem3 and \therosteritem4 together clearly imply 
\therosteritem2. 

\therosteritem2 $\Longrightarrow$ \therosteritem1. 
This is a consequence of the following result.
\newline
{\bf Claim} {\sl Consider two smooth curves $\al,\be:\Bbb R\to G$ with 
$\al(0)=e=\be(0)$ 
which satify the two differential equations
$$\align
\tfrac d{dt}\al(t) &= L_X(\al(t))\\
\tfrac d{dt}\be(t) &= R_X(\be(t)).\\
\endalign$$
Then $\al=\be$ and it is a 1-parameter subgroup.}

We have $\al=\be$ since
$$\align
\tfrac d{dt}(\al(t)\be(-t)) 
&= T\mu^{\be(-t)}.L_X(\al(t)) - T\mu_{\al(t)}.R_X(\be(-t)) \\
&= T\mu^{\be(-t)}.T\mu_{\al(t)}.X - T\mu_{\al(t)}.T\mu^{\be(-t)}.X = 
0.
\endalign
$$
Next we calculate for fixed $s$
$$
\tfrac d{dt}(\be(t-s)\be(s)) 
= T\mu^{\be(s)}.R_X(\be(t-s)) = R_X(\be(t-s)\be(s)).
$$
Hence by the first part of the proof $\be(t-s)\be(s)=\al(t)=\be(t)$.

\smallskip\noindent
The statement about uniqueness follows from \nmb!{2.9}, 
or from the claim.
\qed\enddemo

\subhead \nmb.{3.7}. Definition \endsubhead
Let $G$ be a Lie group with Lie algebra $\g$. We say that $G$ admits 
an \idx{\it exponential mapping} if there exists a smooth mapping
$\exp:\frak g \to G$ such that $t\mapsto \exp(tX)$ is the (unique by 
\nmb!{3.6}) 
1-parameter subgroup with tangent vector $X$ at 0. Then we have by 
\nmb!{3.6}
\roster
\item $\Fl^{L_X}(t,x) = x.\exp(tX)$.
\item $\Fl^{R_X}(t,x) = \exp(tX).x$.
\item $\exp(0)=e$ and $T_0\exp = Id:T_0\frak g=\frak g \to
        T_eG=\frak g$ since $T_0\exp.X = \tfrac d{dt}|_0 \exp(0+t.X) =
        \tfrac d{dt}|_0 \Fl^{L_X}(t,e) = X$.
\item Let $\ph:G\to H$ be a smooth
     homomorphism of between Lie groups admitting exponential 
     mappings. Then the diagram 
$$\newCD 
\frak g   @()\L{\ph'}@(1,0) @()\L{\exp^G}@(0,-1)  
     &  \frak h @()\l{\exp^H}@(0,-1)\\
G @()\L{\ph}@(1,0) &   H  
\endnewCD$$
     commutes, since $t\mapsto \ph(\exp^G(tX))$ is a one parameter
     subgroup of $H$ and $\tfrac d{dt}|_0\ph(\exp^GtX) = \ph'(X)$, so
     $\ph(\exp^GtX)=\exp^H(t\ph'(X))$.
\endroster

\subhead \nmb.{3.8}. Remarks \endsubhead 
\cit!{14}, \cit!{15} gave conditions under which a 
smooth Lie group modelled on Fr\'echet spaces admits exponential 
mappings. We shall elaborate on this notion in \nmb!{5.3} below.
They called this `regular Fr\'echet Lie groups'. We do 
not know of any smooth Fr\'echet Lie group which does not admit an 
exponential mapping. 

If $G$ admits an exponential mapping,
it follows from \nmb!{3.7}.\therosteritem3 that $\exp$ is a 
diffeomorphism from a neighborhood of 0 in $\g$ onto a neighborhood 
of $e$ in $G$, if a suitable inverse function theorem is applicable. 
This is true for example for smooth Banach Lie groups, also for gauge 
groups, but it is wrong for diffeomorphism groups, see 
\cit!{5}. 

If $E$ is a Banach space, then in the Banach Lie group $GL(E)$ 
of all bounded linear automorphisms of $E$ the exponential 
mapping is given by the von Neumann series 
$\exp(X)=\sum_{i=0}^\infty\tfrac1{i!}X^i$. 

If $G$ is connected with exponential 
mapping and $U\subset
\frak g$ is open with $0\in U$, then one may ask whether the group 
generated by $\exp(U)$ equals $G$. 
Note that this is a normal subgroup. So if $G$ is simple, the answer 
is yes.
This is true for connected components of diffeomorphism 
groups and many of their important subgroups.

Results on weakened versions of the Baker-Campbell-Hausdorff formula 
can be found in \cit!{19}.

\subhead \nmb.{3.9}. The adjoint representation \endsubhead 
Let $G$ be a Lie group with Lie algebra $\g$. 
For $a\in G$ we define $\conj_a:G\to G$ by $\conj_a(x)=axa\i$. It
is called the \idx{\it conjugation} or the \idx{\it inner automorphism}
by $a\in G$. This defines a smooth action of $G$ on itself by 
automorphisms.

The adjoint representation $\Ad:G\to GL(\frak g)\subset L(\g,\g)$ 
is given by 
$\Ad(a)=(\conj_a)'=T_e(\conj_a):\frak g\to \frak g$ for $a\in G$. 
By \nmb!{3.5} $\Ad(a)$ is a Lie algebra homomorphism.
By \nmb!{3.2} we have
$\Ad(a)=T_e(\conj_a) = T_a(\mu^{a\i}).T_e(\mu_a) =
     T_{a\i}(\mu_a).T_e(\mu^{a\i})$.

Finally we define the (lower case) \idx{\it adjoint representation} of
the Lie algebra $\frak g$,
$\ad:\frak g\to \frak g\frak l(\frak g):=L(\frak g,\frak g)$,
by $\ad:=\Ad'=T_e\Ad$.

We shall also use the {\it right} \idx{\it Maurer-Cartan form} 
$\ka^r\in\Om^1(G,\g)$, given by 
$\ka^r_g=T_g(\mu^{g\i}):T_gG\to \g$; similarly the 
{\it left Maurer-Cartan form} $\ka^l\in\Om^1(G,\g)$ is 
given by $\ka^l_g=T_g(\mu_{g\i}):T_gG\to \g$.

\proclaim{ Lemma} 
\roster
\item $L_X(a) = R_{\Ad(a)X}(a)$ for $X\in \frak g$
        and $a\in G$.
\item $\ad(X)Y = [X,Y]$ for $X,Y\in \frak g$.
\item $d\Ad = (\ad\o\ka^r).\Ad = 
       \Ad.(\ad\o \ka^l):TG\to L(\g,\g)$.
\endroster
\endproclaim

\demo{Proof} \therosteritem1. $L_X(a)=T_e(\mu_a).X =
T_e(\mu^a).T_e(\mu^{a\i}\o\mu_a).X = R_{\Ad(a)X}(a)$.

\therosteritem2. 
We need some preparations.
Let $V$ be a convenient 
vector space. For $f\in C^\infty(G,V)$ we define the 
\idx{\it left trivialized derivative} 
$D_lf\in C^\infty(G,L(\g,V))$ by 
$$
D_lf(x).X := df(x).T\mu_x.X = (L_Xf)(x).\tag4
$$
For $f\in C^\infty(G,\Bbb R)$ and $g\in C^\infty(G,V)$ we have
$$\align
D_l(f.g)(x).X &= d(f.g)(T_e\mu_x.X)\tag5\\ 
&= df(T_e\mu_x.X).g(x) + f(x).dg(T_e\mu_x.X)\\
&= (f.D_lg + D_lf\otimes g)(x).X.
\endalign$$
 From the fomula
$$\align
D_lD_lf(x)(X)(Y) & = 
     D_l(D_lf(\quad).Y)(x).X \\
&= D_l(L_Yf)(x).X = L_XL_Yf(x).
\endalign$$
follows
$$
D_lD_lf(x)(X)(Y)- D_lD_lf(x)(Y)(X) = L_{[X,Y]}f(x) = D_lf(x).[X,Y].\tag6
$$
We consider now the linear isomorphism 
$L:C^\infty(G,\g)\to \X(G)$
given by $L_f(x)=T_e\mu_x.f(x)$ for $f\in C^\infty(G,\g)$. 
If $h\in C^\infty(G,V)$ we get 
$(L_fh)(x)= D_lh(x).f(x)$.
For $f,g\in C^\infty(G,\g)$ and $h\in C^\infty(G,\Bbb R)$
we get in turn, using \thetag5, generalized to the bilinear pairing 
$L(\g,\Bbb R)\x \g\to \Bbb R$,
$$\align
(L_fL_gh)(x) &= 
D_l(D_lh(\quad).g(\quad))(x).f(x)\\
&= D_lD_lh(x)(f(x))(g(x))
     + D_lh(x).D_lg(x).f(x)\\
([L_f,L_g]h)(x) &= D_l^2h(x).(f(x),g(x)) + 
     D_lh(x).D_lg(x).f(x) -\\
&\quad - D_l^2h(x).(g(x),f(x)) - 
     D_lh(x).D_lf(x).g(x) \\
&= D_lh(x).\Bigl([f(x),g(x)]_{\g} 
     + D_lg(x).f(x) - D_lf(x).g(x)\Bigr)\\
[L_f,L_g] &= L\Bigl([f,g]_{\g} + D_lg.f - D_lf.g\Bigr)\tag7
\endalign$$
Now we are able to prove the second assertion of the lemma. For 
$X,Y\in\g$ we will apply \thetag7 to $f(x)=X$ and $g(x)=\Ad(x\i).Y$. 
We have $L_g=R_Y$ by \therosteritem1, and $[L_f,L_g] = [L_X,R_Y]=0$ by 
\nmb!{3.4}. So
$$\align
0 = [L_X,R_Y](x) &= [L_f,L_g](x) \\
&= L([X,(\Ad\o\nu)Y]_\g + 
D_l((\Ad\o \nu)(\quad).X).Y - 0)(x)\\
[X,Y]& = [X,\Ad(e)Y]=-D_l((\Ad\o \nu)(\quad).X)(e).Y \\
&= d(\Ad(\quad).X)(e).Y = \ad(X)Y.
\endalign$$

\therosteritem3. 
Let $X,Y\in\g$ and $g\in G$ and let $c:\Bbb R\to G$ 
be a smooth curve with $c(0)=e$ and $c'(0)=X$. Then we have
$$\align
(d\Ad(R_X(g))).Y &= \dd t|_0 \Ad(c(t).g).Y 
     = \dd t|_0 \Ad(c(t)).\Ad(g).Y \\
&= \ad(X)\Ad(g)Y = (\ad\o\ka^r)(R_X(g)).\Ad(g).Y,
\endalign$$
and similarly for the second formula.
\qed\enddemo

\subhead \nmb.{3.10} \endsubhead 
Let $r:M\x G\to M$ be a right action, so
$\check r:G\to \operatorname{Diff}(M)$ is a group anti homomorphism. 
We will use the following notation:
$r^a:M\to M$ and $r_x:G\to M$, given by $r_x(a)=r^a(x)=r(x,a)=x.a$.

For any $X\in \frak g$ we define the 
\idx{\it fundamental vector field} $\zeta_X=\ze_X^M\in \X(M)$ by 
$\ze_X(x) = T_e(r_x).X = T_{(x,e)}r.(0_x,X)$.

\proclaim{Lemma} In this situation the following assertions hold:
\roster
\item $\ze:\frak g\to\X(M)$ is a Lie algebra homomorphism.
\item $T_x(r^a).\ze_X(x) = \ze_{Ad(a\i)X}(x.a)$.
\item $0_M\x L_X \in \X(M\x G)$ is $r$-related to 
       $\ze_X\in \X(M)$.\qed
\endroster
\endproclaim

\head\totoc \nmb0{4}. Bundles and connections\endhead

\subhead \nmb.{4.1}. Definition \endsubhead
A \idx{\it principal (fiber) bundle} 
$(P,p,M,G)$ is a smooth mapping $p:P\to M$ such that
there exist an open cover $(U_\al)$ of $M$ and fiber respecting 
diffeomorphisms
$\ph_\al: P|U_\al := p\i(U_\al) \to  U_\al\times G$
with $(\ph_\al\o \ph_\be^{-1})(x,g) = (x,\ph_{\al\be}(x).g)$
for a smooth cocycle of transition functions 
$(\ph_{\al\be}:U_{\al\be}:= U_\al\cap U_\be \to  G)$.
This is called a \idx{\it principal bundle atlas}. 

Each principal bundle admits a unique right action $r: P \times
G \to  P$, called the \idx{\it principal right action}, given
by $\ph_\al(r(\ph_\al^{-1}(x,a),g)) = (x,ag)$. Since left and
right translation on $G$ commute, this is well defined. 
We write $r(u,g) = u.g$ when the meaning is clear. The principal
right action is visibly free and for any $u_x \in P_x$ the
partial mapping $r_{u_x} = r(u_x,\quad): G \to  P_x$ is a diffeomorphism
onto the fiber through $u_x$, whose inverse is denoted by 
$\ta_{u_x}: P_x \to  G$. These inverses together give a smooth
mapping $\ta: P \times_M P \to  G$, whose local expression is
$\ta(\ph_\al^{-1}(x,a),\ph_\al^{-1}(x,b)) = a^{-1}.b$. This mapping
is uniquely determined by the implicit equation 
$r(u_x,\ta(u_x,v_x)) = v_x$, thus we also have
$\ta(u_x.g,u'_x.g') = g\i.\ta(u_x,u'_x).g'$ and $\ta(u_x,u_x) = e$.

\subhead \nmb.{4.2}. Principal connections \endsubhead
 Let $(P,p,M,G)$ be
a principal fiber bundle. Let $VP := (Tp)\i(0_M)\to P$ be the vertical 
bundle.
A (general) connection on $P$ is a smooth fiber projection 
$\Ph: TP \to  VP$, viewed
as a 1-form in $\Om^1(P;VP)\subset \Om^1(P;TP)$, which is called
a \idx{\it principal connection} if it is $G$-equivariant for the
principal right action $r:P\times G \to  P$, so that $T(r^g).\Ph =
\Ph.T(r^g)$  and $\Ph$ is $r^g$-related to itself, or
$(r^g)^*\Ph = \Ph$, for all $g \in G$. 
Then the kernel of $\Ph$ is called the \idx{\it horizontal 
subbundle}, a splitting vector subbundle of $TP\to P$ complementary to 
$VP$.

The vertical bundle of $P$ is trivialized
as a vector bundle over $P$ by the principal action. So 
$\om(X_u) := T_e(r_u)\i.\Ph(X_u) \in \frak g$ is well defined, 
and in this way we
get a $\frak g$-valued 1-form $\om \in \Om^1(P;\frak g)$, which is
called the \idx{\it (Lie algebra valued) connection form} of the
connection $\Ph$. 
Recall from \nmb!{3.10}. the fundamental vector field mapping 
$\ze:\frak g \to  \X(P)$ for the principal right action, which 
trivializes the vertical bundle $P\x \g \cong VP$.
The defining equation for $\om$ can be written also as 
$\Ph(X_u)=\ze_{\om(X_u)}(u)$.

\proclaim{ Lemma} If $\Ph \in \Om^1(P;VP)$ is a principal
connection on the principal fiber bundle $(P,p,M,G)$ then the
connection form has the following two properties:
\roster
\item $\om$ reproduces the generators of fundamental vector
     fields, so that we have $\om(\ze_X(u)) = X$ for all $X \in \frak g$.
\item $\om$ is $G$-equivariant, i.e\. $((r^g)^*\om)(X_u) :=
     \om(T_u(r^g).X_u) = \Ad(g\i).\om(X_u)$ for all $g \in G$ and $X_u
     \in T_uP$. 
\item For the Lie derivative we have $\Cal L_{\ze_X}\om = -\ad(X).\om$.
\endroster 
Conversely a 1-form $\om \in \Om^1(P,\frak g)$ satisfying
\therosteritem1 defines a connection $\Ph$ on $P$ by $\Ph(X_u) =
T_e(r_u).\om(X_u)$, which is a principal connection if and only
if \therosteritem2 is satisfied.
\endproclaim

\demo{Proof} \therosteritem1. $T_e(r_u).\om(\ze_X(u)) =
\Ph(\ze_X(u)) = \ze_X(u) = T_e(r_u).X$. Since $T_e(r_u):\frak g 
\to  V_uP$ is an isomorphism, the result follows. 

\therosteritem2. 
Both directions follow from 
$$\align
T_e(r_{ug}).\om(T_u(r^g).X_u) &=
\ze_{\om(T_u(r^g).X_u)}(ug) = \Ph(T_u(r^g).X_u)\\
T_e(r_{ug}).\Ad(g\i).\om(X_u) &= \ze_{\Ad(g\i).\om(X_u)}(ug) =
T_u(r^g).\ze_{\om(X_u)}(u) \\
&= T_u(r^g).\Ph(X_u).
\endalign$$
\therosteritem3. Let $g(t)$ be a smooth curve in 
$G$ with $g(0)=e$ and $\dd t|_0g(t)=X$. Then $\ph_t:=r^{g(t)}$ is 
a smooth curve of diffeomorphisms on $P$ with $\dd t|_0\ph_t = 
\ze_X$, and by lemma \nmb!{2.10} we have  
$$
\L_{\ze_X}\om = \dd t|_0 (r^{g(t)})^*\om = \dd t|_0 Ad(g(t)\i)\om = 
-ad(X)\om.\qed
$$
\enddemo

\subhead \nmb.{4.3}. Curvature \endsubhead
Let $\Ph$ be a principal
connection on the principal fiber bundle $(P,p,M,G)$ with
connection form $\om\in \Om^1(P;\frak g)$. 

Let us now define the curvature as the obstruction against 
integrability of the horizontal subbundle, i\.e\.
$\Cal R(X,Y):=\Ph[X-\Ph X,Y-\Ph Y]$ for vector fields $X, Y$ on $P$. 
One can check easily that $\Cal R$ is a skew-symmetric bilinear 
$C^\infty(P,\Bbb R)$-module homomorphism, and that 
$(r^{g\i})^*.\Cal R(X,Y)=\Cal R((r^g)^*X,(r^g)^*Y)$ holds, i\.e\. 
$(r^g)^*\Cal R=\Cal R$ for all $g \in G$. Since $\Cal R$ has
vertical values we may again define a $\frak g$-valued 2-form by 
$\Om(X,Y)(u) := - T_e(r_u)\i.\Cal R(X,Y)(u)$, which is called the 
\idx{\it (Lie algebra-valued) curvature form} 
of the connection. 
We also have $\Cal R(X,Y)(u)=-\ze_{\Om(X,Y)(u)}(u)$.
We take the negative sign
here to get in finite dimensions 
the usual curvature form as in \cit!{8}.

We equip the space $\Om(P;\frak g)$ of all $\frak g$-valued forms
on $P$ in a canonical way with the structure of a graded Lie
algebra by 
$$\multline 
[\Ps,\Th]_{\wedge}(\row X1{p+q}) = \\
= \frac 1{p!\,q!} \sum_{\si} \text{sign}\si\,[\Ps(\row
X{\si1}{\si p}),\Th(\row X{\si(p+1)}{\si(p+q)})]_{\frak g}
\endmultline $$
or equivalently by $[\ps\otimes X,\th \otimes Y]_\wedge  := \ps \wedge
\th \otimes [X,Y]_{\frak g}$. 
 From the latter description it is clear that 
$d[\Ps,\Th]_\wedge =[d\Ps,\Th]_\wedge +(-1)^{\deg\Ps}[\Ps,d\Th]_\wedge$.
In particular for $\om \in
\Om^1(P;\frak g)$ we have 
$[\om,\om]_{\wedge}(X,Y) = 2[\om(X),\om(Y)]_{\frak g}$.

\proclaim{ Theorem} The curvature form $\Om$ of a principal
connection with connection form $\om$ has the following
properties: 
\roster
\item $\Om$ is horizontal, i.e. it kills vertical vector fields.
\item The Maurer-Cartan formula holds: $\Om = d\om + \frac 12
     [\om,\om]_{\wedge}\in \Om^2(P;\frak g)$. 
\item $\Om$ is $G$-equivariant in the following sense:
     $(r^g)^*\Om = \Ad(g\i).\Om$. Consequently $\Cal L_{\ze_X}\Om =
     -\ad(X).\Om$. 
\endroster
\endproclaim

\demo{Proof} \therosteritem1 is true for $\Cal R$ by definition. 
For \therosteritem2 we show that the formula holds if at least one vector 
field is vertical, or if both are horizontal. 
For $X \in \frak g$ we have $i_{\ze_X}\Cal R = 0$ by
\therosteritem1, and using \nmb!{4.2}.(\therosteritem1 and 
\therosteritem3) we get
$$\align
i_{\ze_X}(d\om + \frac12[\om,\om]_\wedge ) &= i_{\ze_X}d\om  +
     \frac 12 [i_{\ze_X}\om,\om]_\wedge  - \frac 12 
     [\om,i_{\ze_X}\om]_\wedge  = \\
&= \Cal L_{\ze_X}\om + [X,\om]_\wedge  = - \ad(X)\om + \ad(X)\om = 0.
\endalign $$
So the formula holds for vertical vectors, and for horizontal
vector fields $X,Y$ we have
$$\align 
\Cal R(X,Y) &= \Ph[X-\Ph X,Y-\Ph Y] = \Ph[X,Y] = \ze_{\om([X,Y])} \\ 
(d\om+\frac12 [\om,\om])(X,Y) &= X\om(Y) - Y\om(X) -\om([X,Y]) + 0 =
     -\om([X,Y]). 
\endalign $$
That $\Om$ is really a `tensorial' 2-form follows either from 
\therosteritem2 or from \nmb!{4.4}.\therosteritem4 below.
\qed\enddemo

\subhead \nmb.{4.4}. Local descriptions of principal connections
\endsubhead
We consider a principal fiber bundle $(P,p,M,G)$ with some principal
fiber bundle atlas $(U_\al,\ph_\al:P|U_\al\to U_\al\x G)$ and
corresponding cocycle $(\ph_{\al\be}:U_{\al\be}\to G)$ of
transition functions. 
Let $\Ph = \ze\o\om\in \Om^1(P;VP)$ be a principal connection with 
connection form $\om\in \Om^1(P;\frak g)$. 

We consider the sections $s_\al\in C^\infty(P|U_\al)$ which are
given by $\ph_\al(s_\al(x)) = (x,e)$ and satisfy
$s_\al.\ph_{\al\be}=s_\be$. Then we may associate 
to the connection the collection of the 
$\om_\al:= s_\al{}^*\om\in \Om^1(U_\al;\frak g)$, the
physicists version of the connection. 

\proclaim{ Lemma} These local data have the following properties
and are related by the following formulas.
\roster
\item The forms $\om_\al\in \Om^1(U_\al;\frak g)$
    satisfy the transition formulas
    $$\om_\al = \Ad(\ph_{\be\al}\i)\om_\be + (\ph_{\be\al})\ka^l,$$
    where $\ka^l\in \Om^1(G;\g)$ is the left Maurer-Cartan form from 
    \nmb!{3.9}.
\item The local expression of $\om$ is given by
$$(\ph_\al\i){}^*\om(\xi_x,T\mu_g.X) = (\ph_\al\i){}^*\om(\xi_x,0_g) + X =
    \Ad(g\i)\om_\al(\xi_x) + X.$$
\item The local expression of $\Ph$ is given by 
$$
((\ph_\al)^{-1})^*\Ph(\xi_x,\et_g) 
     = T_e(\mu^g).\om_\al(\xi_x) + \et_g 
     = R_{\om_\al(\xi_x)}(g)  + \et_g
$$
     for $ \xi_x \in T_xU_\al$ and  $\et_g \in T_gG$.
\item The local expression of the curvature $\Cal R$ is given by 
$$
((\ph_\al)\i)^*\Cal R = -R_{d\om_\al+\tfrac12[\om_\al,\om_\al]_\g^\wedge}
$$ 
     so that $\Cal R$ and $\Om$ are indeed `tensorial' 2-forms. 
\endroster
\endproclaim

\demo{Proof} 
For \therosteritem1 to \therosteritem3 plug into the definitions. 
For \therosteritem4 note that the right trivialization or framing 
$(\ka^r,\pi_G):TG\to \g\x G$ 
induces the isomorphism $R:C^\infty(G,\g)\to \X(G)$, given by 
$R_X(x)=T_e(\mu^x).X(x)$.
For the Lie bracket we then have
$$\gather
[R_X,R_Y] = R_{-[X,Y]_\g + dY.R_X - dX.R_Y},\\
R\i[R_X,R_Y] = -[X,Y]_\g + R_X(Y) - R_Y(X).\\
\endgather$$
We write a vector field on $U_\al\x G$ as $(\xi,R_X)$ where 
$\xi:G\to \X(U_\al)$ and $X\in C^\infty(U_\al\x G,\g)$.
Then the local expression of the curvature is given by
$$\align
(\ph_\al&\i)^*\Cal R((\xi,R_X),(\et,R_Y)) =
     (\ph_\al\i)^*(\Cal R((\ph_\al)^*(\xi,R_X),(\ph_\al)^*(\et,R_Y)))\\ 
&= (\ph_\al\i)^*(\Ph[(\ph_\al)^*(\xi,R_X)
     -\Ph(\ph_\al)^*(\xi,R_X),\dots])\\ 
&= (\ph_\al\i)^*(\Ph[(\ph_\al)^*(\xi,R_X)
     -(\ph_\al)^*(R_{\om_\al(\xi)}+R_X),\dots])\\ 
&= (\ph_\al\i)^*(\Ph(\ph_\al)^*[(\xi,-R_{\om_\al(\xi)}),
     (\et,-R_{\om_\al(\et)})])\\ 
&= ((\ph_\al\i)^*\Ph)([\xi,\et]_{\X(U_\al)} -R_{\om_\al(\xi)}(\et) 
     +R_{\om_\al(\et)}(\xi),\\
&\qquad     -\xi(R_{\om_\al(\et)}) + \et(R_{\om_\al(\xi)})
     +R_{-[\om_\al(\xi),\om_\al(\et)] + R_{\om_\al(\xi)}(\om_\al(\et)) 
     -R_{\om_\al(\xi)}(\om_\al(\et))})\\
&= R_{\om_\al([\xi,\et]_{\X(U_\al)} -R_{\om_\al(\xi)}(\et) 
     +R_{\om_\al(\et)}(\xi))} 
     -R_{\xi(\om_\al(\et))} + R_{\et(\om_\al(\xi))}\\
&\quad     +R_{-[\om_\al(\xi),\om_\al(\et)] + R_{\om_\al(\xi)}(\om_\al(\et)) 
     -R_{\om_\al(\xi)}(\om_\al(\et))}\\
&= -R_{(d\om_\al+\tfrac12[\om_\al,\om_\al]_\g^\wedge)(\xi,\et)}.\qed
\endalign$$
\enddemo

\head\totoc \nmb0{5}. Regular Lie groups \endhead

\subhead \nmb.{5.1}. The right and left logarithmic derivatives \endsubhead
Let $M$
be a manifold and let $f:M\to G$ be a smooth mapping into a Lie
group $G$ with Lie algebra $\frak g$. We define the mapping 
$\de^r f:TM\to \frak g$ by the formula
$$
\de^r f(\xi_x) := T_{f(x)}(\mu^{f(x)\i}).T_xf.\xi_x 
\text{ for }\xi_x\in T_xM.
$$
Then $\de^r f$ is a $\frak g$-valued 1-form on $M$, $\de^r f\in
\Om^1(M;\frak g)$. We call $\de^r f$ the
\idx{\it right logarithmic derivative} of $f$, since for $f:\Bbb
R\to(\Bbb R^+,\cdot)$ we have 
$\de^r f(x).1=\frac{f'(x)}{f(x)}=(\log\o f)'(x)$.

Similarly the
\idx{\it left logarithmic derivative} 
$\de^lf\in\Om^1(M,\frak g)$ 
of a smooth mapping $f:M\to G$ is given by 
$$\de^lf.\xi_x= T_{f(x)}(\mu_{f(x)\i}).T_xf.\xi_x.$$

\proclaim{ Lemma} Let $f,g:M\to G$ be smooth. Then the Leibniz rule 
holds:
$$\de^r(f.g)(x) 
     = \de^r f(x) + \Ad(f(x)).\de^r g(x).$$
Moreover the differential form $\de^rf\in\Om^1(M;\g)$ 
satifies the `left Maurer-Cartan equation' (left because it stems 
from the left action of $G$ on itself) 
$$\gather
d\de^rf(\xi,\et) - 
[\de^rf(\xi),\de^rf(\et)]^\g=0,\\
\text{ or }\quad d\de^rf - \frac12 [\de^rf,\de^rf]^\g_\wedge=0,
\endgather$$
where $\xi,\et\in T_x M$, and 
where for $\ph\in\Om^p(M;\g),\ps\in\Om^q(M;\g)$ one puts
$$[\ph,\ps]^\g_\wedge(\xi_1,\dots,\xi_{p+q}) 
     := \frac1{p!q!}\sum_\si\sign(\si)
     [\ph(\xi_{\si1},\dots),\ps_{\si(p+1)},\dots)]^\g.$$
For the left logarithmic derivative 
the corresponding Leibniz rule is uglier, 
and it satisfies the `right Maurer Cartan equation':
$$\gather
\de^l(fg)(x) = \de^lg(x) + 
Ad(g(x)\i)\de^lf(x),\\
d\de^lf + \frac12 
[\de^lf,\de^lf]^\g_\wedge=0.
\endgather$$
\endproclaim

For `regular Lie groups' we will prove a converse to this statement 
later in \nmb!{7.2}.

\demo{Proof}
We treat only the right logarithmic derivative, the proof for the left 
one is similar.
$$\align 
\de^r(f.g)(x) &= T(\mu^{g(x)\i.f(x)\i}).T_x(f.g)  \\
&= T(\mu^{f(x)\i}).T(\mu^{g(x)\i}).T_{(f(x),g(x))}\mu.(T_xf,T_xg)\\
&= T(\mu^{f(x)\i}).T(\mu^{g(x)\i}).\left(T(\mu^{g(x)}).T_xf +
     T(\mu_{f(x)}).T_xg\right) \\
&= \de^r f(x) + \Ad(f(x)).\de^r g(x).
\endalign$$

We shall use now principal bundle geometry from section \nmb!{3}.
We consider the trivial principal bundle 
$\operatorname{pr_1}:M\x G\to M$ with right principal action. Then 
the submanifolds $\{(x,f(x).g):x\in M\}$ for $g\in G$ form a 
foliation of $M\x G$ whose tangent distribution is complementary to the 
vertical bundle $M\x TG \subseteq T(M\x G)$ and is invariant under the 
principal right $G$-action. So it is the horizontal distribution of a 
principal connection on $M\x G\to G$. For a tangent vector 
$(\xi_x,Y_g)\in T_xM\x T_gG$ the horizontal part is the right translate 
to the foot point $(x,g)$ of $(\xi_x,T_xf.\xi_x)$, so the decomposition in 
horizontal and vertical parts according to this distribution is
$$
(\xi_x,Y_g) = (\xi_x, T(\mu^g).T(\mu^{f(x)\i}).T_xf.\xi_x)
+(0_x,Y_g-T(\mu^g).T(\mu^{f(x)\i}).T_xf.\xi_x).
$$
Since the fundamental vector fields for the right action on $G$ are 
the left invariant vector fields, the corresponding connection form 
is given by
$$\align
\om^r(\xi_x,Y_g) &= 
T(\mu_{g\i}).(Y_g-T(\mu^g).T(\mu^{f(x)\i}).T_xf.\xi_x),\\
\om^r_{(x,g)} &= T(\mu_{g\i}) - \Ad(g\i).\de^r f_x,\\
\om^r &= \ka^l - (\Ad\o\nu).\de^r f,\tag1
\endalign$$
where $\ka^l:TG\to\g$ is the left Maurer-Cartan form on $G$ (the 
left trivialization), given by $\ka^l_g=T(\mu_{g\i})$. Note that 
$\ka^l$ is the principal connection form for the (unique) principal 
connection $p:G\to \text{point}$ with right principal action, which 
is flat so that the right (from right action) Maurer-Cartan equation 
equation holds in the form
$$d\ka^l+\tfrac12[\ka^l,\ka^l]_\wedge=0.\tag2$$

The principal connection $\om^r$ is flat since we got it via the horizontal 
leaves, so the principal connection form vanishes:
$$\align
0&= d\om^r+\tfrac12[\om^r,\om^r]_\wedge\tag3\\ 
&= d\ka^l + \tfrac12[\ka^l,\ka^l]_\wedge 
     - d(\Ad\o\nu) \wedge \de^r f - 
     (\Ad\o\nu). d\de^r f \\
&\quad - [\ka^l,(\Ad\o\nu). \de^r f]_\wedge + 
     \tfrac12[(\Ad\o\nu).\de^r f,
     (\Ad\o\nu). \de^r f]_\wedge\\
&= - (\Ad\o\nu). (d\de^r f - 
     \tfrac12[\de^r f,\de^r f]_\wedge), 
\endalign$$ 
where we used \thetag2
and since for $\xi\in\g$ and a smooth curve $c:\Bbb R\to G$ with 
$c(0)=e$ and $c'(0)=\xi$ we have:
$$\align
d(\Ad\o\nu) (T(\mu_g)\xi) &= \ddt \Ad(c(t)\i.g\i) 
     = -\ad(\xi)\Ad(g\i) \\ 
&= -\ad\Bigl(\ka^l(T(\mu_g)\xi)\Bigr)(\Ad\o\nu)(g),\\
d(\Ad\o\nu) 
&=-(\ad\o\ka^l).(\Ad\o\nu).\tag4 
\endalign$$
So we have $d\de^r f - \tfrac12[\de^r f,\de^r f]_\wedge$ as asserted.

For the left logarithmic derivative $\de^l f$ the proof is 
similar, and we discuss only the essential deviations. First note 
that on the trivial principal bundle $\operatorname{pr_1}:M\x G\to M$ 
with left principal action of $G$ the fundamental vector fields are 
the right invariant vector fields on $G$, and that for a principal 
connection form $\om^l$ the curvature form is given by 
$d\om^l-\frac12[\om^l,\om^l]_\wedge$. Look at the proof of 
theorem \nmb!{4.3}  to see this. 
The connection form is then given by
$$
\om^l = \ka^r - \Ad.\de^l f,\tag1'
$$
where the right Maurer-Cartan form 
$(\ka^r)_g=T(\mu^{g\i}):T_gG\to\g$ now satifies the left 
Maurer-Cartan equation
$$d\ka^r-\frac12[\ka^r,\ka^r]_\wedge = 0. \tag2'$$
Flatness of $\om^l$ now leads to the computation
$$\align
0&= d\om^l-\tfrac12[\om^l,\om^l]_\wedge\tag3'\\ 
&= d\ka^r 
     - \tfrac12[\ka^r,\ka^r]_\wedge      
     - d\Ad \wedge \de^l f - \Ad. d\de^l f \\
&\quad + [\ka^r,\Ad. \de^l f]_\wedge - 
     \tfrac12[\Ad.\de^l f,\Ad. \de^l f]_\wedge\\
&= - \Ad. (d\de^l f 
     + \tfrac12[\de^l f,\de^l f]_\wedge), 
\endalign$$ 
where we used $d\Ad =(\ad\o\ka^r)\Ad$ from \nmb!{3.9}.(3) 
directly.
\qed\enddemo

\subhead\nmb.{5.2} \endsubhead
Let $G$ be a Lie group with Lie algebra $\g$. 
For a closed interval $I\subset \Bbb R$  and for 
$X\in C^\infty(I,\g)$ we consider the ordinary differential equation
$$
\cases g(t_0)=e \\
       \dd tg(t)=T_e(\mu^{g(t)})X(t) = R_{X(t)}(g(t)),
          \quad\text{ or }\ka^r(\dd tg(t)) = X(t),
\endcases\tag1
$$
for local smooth curves $g$ in $G$, where $t_0\in I$.

\proclaim{Lemma}  
\roster
\item[2] Local solution curves $g$ of the differential equation 
       \thetag1 are uniquely determined.  
\item  If for fixed $X$ the differential equation \thetag1 has a 
       local solution near each $t_0\in I$, then it has also a global 
       solution $g\in C^\infty(I,G)$.  
\item  If for all $X\in C^\infty(I,\g)$ the differential equation 
       \thetag1 has a local solution near one fixed $t_0\in I$, then 
       it has also a global solution $g\in C^\infty(I,G)$ for each 
       $X$. Moreover, if the local solutions near $t_0$ depend 
       smoothly on the vector fields $X$ (see the proof for the exact 
       formulation), then so does the global solution. 
\item  The curve $t\mapsto g(t)\i$ is the unique local smooth curve $h$ in 
       $G$ which satifies  
$$
\cases h(t_0)=e \\
\dd th(t) = T_e(\mu_{h(t)})(-X(t)) = L_{-X(t)}(h(t)),
          \quad\text{ or }\ka^l(\dd th(t)) = -X(t).
\endcases$$
\endroster
\endproclaim

\demo{Proof}
\therosteritem2.
Suppose that $g(t)$ and $g_1(t)$ both satisfy \thetag1. Then on 
the intersection of their intervals of definition we have 
$$\align
\dd t(g(t)\i\,g_1(t)) &= 
-T(\mu^{g_1(t)}).T(\mu_{g(t)\i}).T(\mu^{g(t)\i}).T(\mu^{g(t)}).X(t)\\
&\quad + T(\mu_{g(t)\i}).T(\mu^{g_1(t)}).X(t) = 0,
\endalign$$
so that $g=g_1$. 

\therosteritem3. 
It suffices to prove the claim for every compact subintervall of 
$I$, so let $I$ be compact.
If $g$ is a local solution of \thetag1 then $t\mapsto g(t).x$ is a 
local solution of the same differential equation with initial value 
$x$. By assumption for each $s\in I$ there is 
a unique solution $g_s$ of the differential 
equation with $g_s(s)=e$; so there 
exists $\de_s>0$ such that
$g_s(s+t)$ is defined for $|t|<\de_s$. 
Since $I$ is compact there exist 
$s_0<s_1<\dots<s_k$ such that $I=[s_0,s_k]$ and $s_{i+1}-s_i<\de_{s_i}$. 
Then we put
$$
g(t):=\cases g_{s_0}(t) \quad & \text{ for }s_0\le t \le s_1\\
             g_{s_1}(t).g_{s_0}(s_1) \quad & \text{ for }s_1\le t \le s_2\\       
               \dotsc & \\
             g_{s_i}(t).g_{s_{i-1}}(s_i)\dots g_{s_0}(s_1) \quad 
                    & \text{ for }s_i\le t \le s_{i+1}\\
               \dotsc & 
     \endcases
$$
which is smooth by the first case and solves the problem.

\therosteritem4. 
Given $X:I\to \g$ we first extend $X$ to a smooth curve 
$\Bbb R\to \g$, 
using the formula of \cit!{18}.
For $t_1\in I$, by assumption there exists a local 
solution $g$ near $t_0$ of the translated vector field
$t\mapsto X(t_1-t_0+t)$, thus $t\mapsto g(t_0-t_1+t)$ is a solution 
near $t_1$ of $X$. So by assertion \therosteritem3 the differential 
equation has a global solution for $X$ on $I$.

Now we assume that the local solutions near $t_0$ depend smoothly in 
the vector field: So for any smooth curve 
$X:\Bbb R\to C^\infty(I,\g)$ we have:
\block
For each compact intervall $K\subset \Bbb R$ there is a neighborhood 
$U_{X,K}$ of $t_0$ in $I$ and a smooth mapping $g:K\x U_{X,K}\to G$ 
with
$$
\cases g(k,t_0)=e \\
       \dd tg(k,t)=T_e(\mu^{g(k,t)}).X(k)(t)\quad\text{ for all 
       }k\in K, t\in U_{X,K}. 
\endcases
$$
\endblock
Given a smooth curve $X:\Bbb R\to C^\infty(I,\g)$ we extend (or lift) it 
smoothly to $X:\Bbb R\to C^\infty(\Bbb R,\g)$ by 
using the formula of \cit!{18}.
Then the smooth parameter $k$ 
from the compact intervall $K$ passes smoothly through the 
proofs given above to give a smooth global solution $g:K\x I\to G$.
So the `solving operation' respects smooth curves and thus is 
`smooth'.

\therosteritem5. 
One can show in a similar way that $h$ is the unique 
solution of \thetag5 by differentiating $h_1(t).h(t)\i$.
Moreover the curve $t\mapsto g(t)\i=h(t)$ satisfies 
\thetag5, since
$$
\dd t(g(t)\i) = - T(\mu_{g(t)\i}).T(\mu^{g(t)\i}).T(\mu^{g(t)}).X(t) 
= T(\mu_{g(t)\i}).(-X(t)). \qed
$$
\enddemo

\subhead \nmb.{5.3}. Definition. Regular Lie groups \endsubhead
If for each $X\in C^\infty(\Bbb R,\g)$ there exists 
$g\in C^\infty(\Bbb R,G)$ satisfying 
$$
\cases g(0)=e \\
       \dd tg(t)=T_e(\mu^{g(t)})X(t) = R_{X(t)}(g(t)),\\
          \quad\text{ or }\ka^r(\dd tg(t))= \de^rg(\partial_t) = X(t).
\endcases\tag1
$$
then we write 
$$\gather
\evol^r_G(X) = \evol_G(X):=g(1),\\
\Evol^r_G(X)(t) := \evol_G(s\mapsto tX(ts)) = g(t),
\endgather$$ 
and call it the \idx{\it right evolution} 
of the curve $X$ in $G$. By lemma \nmb!{5.2} the 
solution of the differential equation \thetag1 is unique, 
and for global existence it is sufficient that it has a local solution.
Then 
$$
\Evol^r_G: C^\infty(\Bbb R,\g) \to \{g\in C^\infty(\Bbb R,G):g(0)=e\}
$$
is bijective with inverse the right logarithmic derivative $\de^r$.

The Lie group $G$ is called a
\idx{\it regular Lie group} if 
$\evol^r:C^\infty(\Bbb R,\g)\to G$ exists and is smooth.

We also write 
$$\gather
\evol^l_G(X) = \evol_G(X):=h(1),\\
\Evol^l_G(X)(t) := \evol^l_G(s\mapsto tX(ts)) = h(t),
\endgather$$ 
if $h$ is the (unique) solution of 
$$
\cases h(0)=e \\
     \dd th(t) = T_e(\mu_{h(t)})(X(t)) = L_{X(t)}(h(t)),\\
          \quad\text{ or }\ka^l(\dd th(t))=\de^lh(\partial_t) = X(t).
\endcases\tag2$$
Clearly $\evol^l:C^\infty(\Bbb R,\g)\to G$ exists and is 
also smooth if $\evol^r$ does, 
since we have $\evol^l(X)=\evol^r(-X)\i$ by lemma \nmb!{5.2}. 

Let us collect some easily seen properties of the evolution mappings.
If $f\in C^\infty(\Bbb R,\Bbb R)$, then we have 
$$\align
\Evol^r(X)(f(t)) &= \Evol^r(f'.(X\o f))(t).\Evol^r(X)(f(0)),\\
\Evol^l(X)(f(t)) &= \Evol^l(X)(f(0)).\Evol^l(f'.(X\o f))(t).
\endalign$$
If $\ph:G\to H$ is a smooth homomorphism between regular Lie groups 
then the diagram
$$\newCD 
C^\infty(\Bbb R,\frak g)   @()\L{\ph'_*}@(1,0) @()\L{\evol_G}@(0,-1)  &  
     C^\infty(\Bbb R,\frak h) @()\l{\evol_H}@(0,-1)\\ 
G @()\L{\ph}@(1,0) &   H  
\endnewCD$$
commutes, since 
$\dd t\ph(g(t))=T\ph.T(\mu^{g(t)}).X(t)=T(\mu^{\ph(g(t))}).\ph'.X(t)$.

Note that each regular Lie group admits an exponential mapping, 
namely the restriction of $\evol^r$ to the 
constant curves $\Bbb R\to \g$. A Lie group is regular if 
and only if its universal covering group is regular.

This notion of regularity is a weakening of the same notion of 
\cit!{14}, \cit!{15}, who considered a sort of 
product integration property on a  
smooth Lie group modelled on Fr\'echet spaces. Our notion here is 
due to \cit!{13}. Up to now the following statement holds:
\roster
\item "" All known Lie groups are regular.
\endroster
Any Banach Lie group is regular since we may consider the time 
dependent right invariant vector field $R_{X(t)}$ on $G$ and its 
integral curve $g(t)$ starting at $e$, which exists and depends 
smoothly on (a further parameter in) $X$. In particular 
finite dimensional Lie groups are regular.  

For diffeomorphism groups the evolution operator 
is just integration of time dependent 
vector fields with compact support.

\subhead\nmb.{5.4}. Some abelian regular Lie groups \endsubhead
For $(E,+)$, where $E$ is a convenient vector space, we have 
$\evol(X)=\int_0^1X(t)dt$, so convenient vector spaces are regular 
abelian Lie groups. We shall need `discrete' subgroups, which is not 
an obvious notion since $(E,+)$ is not a topological group:
the addition is continuous only $c^\infty(E\x E)\to c^\infty E$, and 
not for the cartesian product of the $c^\infty$-topologies.

Next let $Z$ be a `discrete' subgroup of a convenient vector space 
$E$ in the sense that there exists a $c^\infty$-open neighborhood $U$ 
of zero in $E$ such that $U\cap (z+U)=\emptyset$ for all 
$0\ne z\in Z$ (equivalently $(U-U)\cap (Z\setminus 0)=\emptyset$). For 
that it suffices e\.g\. that $Z$ is discrete in the bornological 
topology on $E$. Then $E/Z$ is an abelian but possibly non Hausdorff 
Lie group. It does not suffice to take $Z$ discrete in the 
$c^\infty$-topology: Take as $Z$ the subgroup generated by $A$ in 
$\Bbb R^{\Bbb N\x c_0}$ in the proof of \cit!{4}, 6.2.8.(iv).

Let us assume that $Z$ fulfills the stronger condition: there exists a 
symmetric $c^\infty$-open neighborhood $W$ of 0 such that 
$(W+W)\cap (z+W+W)=\emptyset$ for all $0\ne z\in Z$
(equivalently $(W+W+W+W)\cap (Z\setminus 0)=\emptyset$).
Then $E/Z$ is Hausdorff and thus an abelian regular Lie group, since 
its universal cover $E$ is regular. Namely, for $x\notin Z$, we have 
to find neighborhoods $U$ and $V$ of 0 such that $(Z+U)\cap (x+Z+V) = 
\emptyset$. There are two cases. If $x\in Z+W+W$ then there is a 
unique $z\in Z$ with $x\in z+W+W$ and we may choose $U,V\subset W$ such 
that $(z+U)\cap (x+V)=\emptyset$; then $(Z+U)\cap (x+Z+V)=\emptyset$.
In the other case, if $x\notin Z+W+W$, then we have 
$(Z+W)\cap (x+Z+W)=\emptyset$.

Notice that the two conditions above and their consequences also hold 
for general (non-abelian) (regular) Lie groups instead of $E$, and 
their `discrete' normal subgroups (which turn out to be central if 
$G$ is connected).

It would be nice if any regular abelian Lie group would be of the 
form $E/Z$ described above. A first result in this direction is that
for an abelian Lie group $G$ with Lie algebra $\g$ which admits a 
smooth exponential mapping $\exp:\g\to G$ one can check easily by 
using \nmb!{5.10} that 
$\dd t(\exp(-tX).\exp(tX+Y)) = 0$ so that $\exp$ is a smooth 
homomorphism of Lie groups.

Let us consider some examples. For the first one we consider a 
discrete subgroup $Z\subset \Bbb R^{\Bbb N}$. There exists a 
neighborhood of 0, without loss of the form 
$U\x \Bbb R^{\Bbb N\setminus n}$ for $U\subset \Bbb R^n$, 
with $U\cap (Z\setminus 0)=\emptyset$. Then we consider the following
diagram of Lie group homomorphisms 
$$\cgaps{1;1;0.4}
\newCD
0 @(1,0) @(0,-1) & \Bbb R^{\Bbb N\setminus n} @(1,0) @(0,-1) &
     \Bbb R^{\Bbb N\setminus n} @(0,-1) & \\
Z @(1,0) @()\L{\cong}@(0,-1) & \Bbb R^{\Bbb N} @(1,0) @()\L{\pi}@(0,-1) &
     \Bbb R^{\Bbb N}/Z @(0,-1) @()\a=@(1,0) & 
     (S^1)^k\x \Bbb R^{\Bbb N\setminus (n-k)} @(0,-1)\\
\pi(Z) @(1,0) & \Bbb R^{n} @(1,0) &
     \Bbb R^{n}/\pi(Z) @()\a=@(1,0) & (S^1)^k\x \Bbb R^{n-k}\\
\endnewCD$$
which has exact lines and columns. For the right hand column we use a 
diagram chase to see this. Choose a global linear section of $\pi$ 
inverting $\pi|Z$. This factors to a global homomorphism of the right 
hand side column.

As next example we consider 
$\Bbb Z^{(\Bbb N)}\subset \Bbb R^{(\Bbb N)}$. Then obviously 
$\Bbb R^{(\Bbb N)}/\Bbb Z^{(\Bbb N)}=(S^1)^{\Bbb N}$, which is 
a smooth (even real analytic, see \cit!{10}) manifold modeled on 
$\Bbb R^{(\Bbb N)}$.
The reader may 
convince himself that the general Lie group covered by 
$\Bbb R^{(\Bbb N)}$ is isomorphic to 
$(S^1)^{(A)}\x \Bbb R^{(\Bbb N\setminus A)}$ 
for $A\subseteq \Bbb N$. 

As another example one may check easily that 
$\ell^\infty/(\Bbb Z^{\Bbb N}\cap \ell^\infty) = (S^1)^{\Bbb N}$, 
equipped with the `uniform box topology'.

\subhead\nmb.{5.5}. Extensions of Lie groups \endsubhead
Let $H$ and $K$ be Lie groups.
A Lie group $G$ is called an smooth \idx{\it extension \ign{  of groups}} of 
$H$ with kernel $K$ if we have a short exact sequence of groups
$$
\{e\} \to K \East{i}{} G \East{p}{} H \to \{e\},\tag1
$$
such that $i$ and $p$ are smooth
and one of the following two equivalent conditions is satified:
\roster
\item[2] $p$ admits a local smooth section $s$ near $e$ (equivalently near 
       any point), and $i$ is initial (i\.e\. any $f$ into $K$ is 
       smooth if and only if $i\o f$ is smooth).
\item $i$ admits a local smooth retraction $r$ near $e$ (equivalently 
       near any point), and $p$ is final (i\.e\. $f$ from $H$ is smooth 
       if and only if $f\o p$ is smooth).
\endroster
Of course by $s(p(x))i(r(x))=x$ the two conditions are equivalent, 
and then $G$ is locally diffeomorphic to $K\x H$ via $(r,p)$ with 
local inverse $(i\o\pr_1).(s\o\pr_2)$.

Not every smooth exact sequence of Lie groups admits local sections 
as required in \therosteritem2. Let for example $K$ be a closed 
linear subspace in a convenient vector space $G$ which is not a 
direct summand, and let $H$ be $G/K$. Then the tangent mapping at 0 
of a local smooth splitting would make $K$ a direct summand.

\proclaim{Theorem} Let $\{e\} \to K \East{i}{} G \East{p}{} H 
\to \{e\}$ be a smooth extension of Lie groups.
Then $G$ is regular if and only if both $K$ and $H$ are regular. 
\endproclaim

\demo{Proof}
Clearly the induced sequence of Lie algebras is also exact,
$$
0 \to \frak k \East{i'}{} \frak g \East{p'}{} \frak h \to 0,
$$
with a bounded linear section $T_es$ of $p'$, so $\frak g$ is 
isomorphic to $\frak k\x\frak h$ as convenient vector space. 

Let us suppose that $K$ and $H$ are regular.
Given $X\in C^\infty(\Bbb R,\g)$, we consider 
$Y(t):=p'(X(t))\in\frak h$ with evolution curve $h$ satisfying
$\dd th(t)=T(\mu^{h(t)}).Y(t)$ and $h(0)=e$. 
By lemma 
\nmb!{5.2} it suffices to find smooth local solutions $g$ near $0$ 
of $\dd t g(t)= T(\mu^{g(t)}).X(t)$ with $g(0)=e$,
depending smoothly on $X$. 
We look for solutions of the form $g(t)=s(h(t)).i(k(t))$, where 
$k$ is a local evolution curve in $K$ of a suitable curve 
$t\mapsto Z(t)$ in $\frak k$, i\.e\.  
$\dd t k(t)=T(\mu^{k(t)}).Z(t)$ and $k(0)=e$. For this ansatz we have 
$$\align
\dd tg(t) &= \dd t\Bigl(s(h(t)).i(k(t))\Bigr) 
     = T(\mu_{s(h(t))}).Ti.\dd tk(t) + T(\mu^{i(k(t))}).Ts.\dd th(t) \\
&=  T(\mu_{s(h(t))}).Ti.T(\mu^{k(t)}).Z(t)
     + T(\mu^{i(k(t))}).Ts.T(\mu^{h(t)}).Y(t), 
\endalign$$
and we want this to be 
$$
T(\mu^{g(t)}).X(t) = T(\mu^{s(h(t)).i(k(t))}).X(t) = 
     T(\mu^{i(k(t))}).T(\mu^{s(h(t))}).X(t).
$$
Using $i\o\mu^{k}=\mu^{i(k)}\o i$ one quickly sees that 
$$
i'.Z(t) := \Ad\Bigl(s(h(t))\i\Bigr).\left(X(t) - 
     T(\mu^{s(h(t))\i}).Ts.T(\mu^{h(t)}).Y(t)  \right) \in \ker p'
$$
solves the problem, so $G$ is regular.

Let now $G$ be regular. If $Y\in C^\infty(\Bbb R,\frak h)$, then 
$p\o\Evol_G^r(s'\o Y)=\Evol_H(Y)$, since for 
$g:=\Evol_G^r(s'\o Y)$ we have
$$\dd t p(g(t))=Tp.\dd tg(t)
     =Tp.T(\mu^{g(t)}).T_es.Y(t)=T(\mu^{p(g(t))}).Y(t).$$ 
If $U\in C^\infty(\Bbb R,\frak k)$, then $p\o\Evol_G(i'\o U) = 
\Evol_H(0)=e$ so that $\Evol_G(i'\o U)(t)\in i(K)$ for all $t$ and 
thus equals $i(\Evol_K(U)(t))$.
\qed\enddemo

\subhead\nmb.{5.6}. Subgroups of regular Lie groups \endsubhead
Let $G$ and $K$ be Lie groups, let $G$ be regular and let $i:K\to G$ 
be a smooth homomorphism which is initial (see \nmb!{5.5}) with
$T_ei=i':\frak k\to \frak g$ injective. 
We suspect that $K$ is then regular, but we are only able to prove 
this under the following assumption.
\block
There is an open neighborhood $U\subset G$ of $e$ and a smooth 
mapping $p:U\to E$ into a convenient vector space $E$ such that 
$p\i(0)=K\cap U$ and $p$ constant on left cosets $Kg\cap U$.
\endblock

\demo{Proof}
For $Z\in C^\infty(\Bbb R,\frak k)$ we consider 
$g(t)=\Evol_G(i'\o Z)(t)\in G$. Then we have 
$\dd t(p(g(t)))=Tp.T(\mu^{g(t)}).i'(Z(t))=0$ 
by the assumption, so $p(g(t))$ is constant $p(e)=0$, thus 
$g(t)=i(h(t))$ for a smooth curve $h$ in $H$, since $i$ is initial.
Then $h=\Evol_H(Y)$ since $i$ is an immersion, and $h$ depends 
smoothly on $Z$ since $i$ is initial.
\qed\enddemo

\subhead\nmb.{5.7}. Abelian and central extensions \endsubhead
 From theorem \nmb!{5.5} it is clear that any smooth extension $G$ of 
a regular Lie group $H$ with an abelian regular Lie group $(K,+)$ is 
again regular. We shall describe $\Evol_G$ in terms of $\Evol_G$, 
$\Evol_K$, and in terms of the action of $H$ on $K$ and the cocycle 
$c:H\x H\to K$ if the latter exists. 

Let us first recall these notions.
If we have a smooth extension with abelian normal subgroup $K$,
$$
\{e\} \to K \East{i}{} G \East{p}{} H \to \{e\},
$$
then a unique smooth action $\al:H\x K\to K$ by automorphisms is given by 
$i(\al_h(k))=s(h)i(k)s(h)\i$, where $s$ is any smooth local section of 
$p$ defined near $h$. If moreover $p$ admits a global smooth section 
$s:H\to G$, which we assume without loss to satisfy $s(e)=e$, 
then we consider the smooth mapping $c:H\x H\to K$ given 
by $ic(h_1,h_2):= s(h_1).s(h_2).s(h_1.h_2)\i$. Via the diffeomorphism
$K\x H\to G$ given by $(k,h)\mapsto i(k).s(h)$ the identity 
corresponds to $(0,e)$, the multiplication and 
the inverse in 
$G$ look as follows:
$$\cases
(k_1,h_1).(k_2,h_2)=(k_1+\al_{h_1}k_2+c(h_1,h_2),h_1h_2),\\
     (k,h)\i=(-\al_{h\i}(k)-c(h\i,h),h\i).
\endcases\tag1$$
Associativity and $(0,e)^2=(0,e)$ correspond to the fact that $c$ 
satisfies the following \idx{\it cocycle condition} and normalization
$$\cases
\al_{h_1}(c(h_2,h_3))-c(h_1h_2,h_3)+c(h_1,h_2h_3)-c(h_1,h_2)=0\\
c(e,e)=0.
\endcases\tag2$$
These imply that $c(e,h)=0=c(h,e)$ and $\al_h(c(h\i,h))=c(h,h\i)$.
For a central extension the action is trivial, $\al_h=\Id_K$ for all 
$h\in H$.

If conversely $H$ acts smoothly by automorphisms on an abelian Lie 
group $K$ and if $c:H\x H\to K$ satisfies \thetag2, then \thetag1 
describes a smooth Lie group structure on $K\x H$, 
which is a smooth extension of $H$ over $K$ with a global smooth 
section.

For later purposes let us compute 
$$\align
(0,h_1).(0,h_2)\i &= (-\al_{h_1}(c(h_2\i,h_2))+ 
     c(h_1,h_2\i),h_1h_2\i),\\
T_{(0,h_1)}(\mu^{(0,h_2)\i}).(0,Y_{h_1}) &= (- T(\al^{c(h_2\i,h_2)}).Y_{h_1} 
     + T(c(\quad,h_2\i)).Y_{h_1}, T(\mu^{h_2\i}).Y_{h_1}).
\endalign$$
Let us now assume that $K$ and $H$ are moreover regular Lie groups.
We consider a curve $t\mapsto X(t)=(U(t),Y(t))$ in the Lie algebra 
$\frak g$ which as convenient vector space  equals 
$\frak k\times\frak h$. 
 From the proof of \nmb!{5.5} we get 
that
$$\align
g(t):&=\Evol_{G}(U,Y)(t) = (0,h(t)).(k(t),e) = 
     (\al_{h(t)}(k(t)),h(t)), \text{ where}\\
h(t):&=\Evol_H(Y)(t) \in H,\\
(Z(t),0):&= 
\Ad_{G}(0,h(t))\i\Bigl((U(t),Y(t))-T\mu^{(0,h(t))\i}.
     (0,\dd t h(t))\Bigr) \\
Z(t) &= T_0(\al_{h(t)\i}).\Bigl(U(t) 
     + \bigl(T(\al^{c(h(t)\i,h(t))})
     - T(c(\quad,h(t)\i))\bigr).\dd t h(t)\Bigr),\\
k(t):&= \Evol_K(Z)(t)\in K.\\
\endalign$$

\subhead{\nmb.{5.8}. Semidirect products}\endsubhead 
 From theorem \nmb!{5.5} we see immediately that the semidirect product 
of regular Lie groups is again regular. Since we shall need explicit 
formulas later we specialize the proof of \nmb!{5.5} to this case.

Let $H$ and $K$ be regular Lie groups with Lie algebras $\frak h$ and 
$\frak k$, respectively. Let $\al:H\x K\to K$ be smooth such that 
$\check\al:H\to \operatorname{Aut}(K)$ is a group homomorphism. Then 
the semidirect product $K\rtimes H$ is the Lie group $K\x H$ with 
multiplication $(k,h).(k',h')=(k.\al_h(k'),h.h')$ and inverse 
$(k,h)\i=(\al_{h\i}(k)\i,h\i)$. We have then
$T_{(e,e)}(\mu^{(k',h')}).(U,Y)
     =(T(\mu^{k'}).U+T(\al^{k'}).Y,T(\mu^{h'}).Y)$.

Now we consider a curve $t\mapsto X(t)=(U(t),Y(t))$ in the Lie algebra 
$\frak k\rtimes\frak h$. Since $s:h\mapsto (e,h)$ is a smooth 
homomorphism of Lie groups, from the proof of \nmb!{5.5} we get 
that
$$\align
g(t):&=\Evol_{K\rtimes H}(U,Y)(t) = (e,h(t)).(k(t),e) = 
(\al_{h(t)}(k(t)),h(t)), \text{ where}\\
h(t):&=\Evol_H(Y)(t) \in H,\\
(Z(t),0):&= \Ad_{K\rtimes H}(e,h(t)\i)(U(t),0) =
(T_e(\al_{h(t)\i}).U(t),0),\\
k(t):&= \Evol_K(Z)(t)\in K.\\
\endalign$$


\proclaim{\nmb.{5.9}. Corollary} Let $G$ be a Lie group.
Then via right trivialization $(\ka^r,\pi_G):TG\to \g\x G$ 
the tangent group $TG$ is isomorphic to the semidirect product 
$\g\rtimes G$, where $G$ acts by $\Ad:G\to\operatorname{Aut}(\g)$.

So if $G$ is a regular Lie group, then $TG\cong \g\rtimes G$ is also 
regular, and $T\evol^r_G$ corresponds to $\evol^r_{TG}$.
In particular for 
$(Y,X)\in C^\infty(\Bbb R,\g\x \g)=TC^\infty(\Bbb R,\g)$, where $X$ 
is the footpoint, and we have 
$$\gather
\evol^r_{\g\rtimes G}(Y,X) = 
\Bigl(\Ad(\evol_G^r(X))\int_0^1\Ad(\Evol^r_G(X)(s)\i).Y(s)\,ds,\; 
     \evol^r_G(X)\Bigr)\\
T_X\evol^r_G.Y = 
T(\mu_{\evol^r_G(X)}).\int_0^1\Ad(\Evol^r_G(X)(s)\i).Y(s)\,ds,\\
T_X(\Evol^r_G(\quad)(t)).Y = T(\mu_{\Evol^r_G(X)(t)}).
     \int_0^t\Ad(\Evol^r_G(X)(s)\i).Y(s)\,ds,
\endgather$$
\endproclaim

Note that in the semidirect product 
representation $TG\cong \g\rtimes G$ the footpoint appears in the 
right factor $G$, contrary to the usual convention. We followed this 
also in $T\g = \g\rtimes \g$.

\demo{Proof}
Via right trivialization the tangent group 
$TG$ is the semidirect product $\g\rtimes G$,
where $G$ acts on the Lie algebra $\g$ by 
$\Ad:G\to \operatorname{Aut}(\g)$, 
because by \nmb!{3.2} we have for $g, h\in G$ and $X,Y\in\g$, where 
$\mu=\mu_G$ is the multiplication on $G$:
$$\align
T_{(g,h)}\mu.(R_X(g),R_Y(h)) &= T(\mu^h).R_X(g) + T(\mu_g).R_Y(h) \\
&= T(\mu^h).T(\mu^g).X + T(\mu_g).T(\mu^h).Y \\
&= R_X(gh) + R_{\Ad(g)Y}(h), \\
T_{g}\nu.R_X(g) &= - T(\mu^{g\i}).T(\mu_{g\i}).T(\mu^g).X \\
&= - R_{\Ad(g\i)X}(g\i),
\endalign$$
so that we have
$$\align
\mu_{\g\rtimes G}((X,g),(Y,h)) &= (X + \Ad(g)Y,gh)\tag1\\
\nu_{\g\rtimes G}(X,g) &= (- \Ad(g\i)X,g\i).
\endalign$$
Now we shall prove that the following diagram commutes and that the 
equations of the corollary follow. The lower triangle commutes by 
definition.
$$\cgaps{;}\newCD
TC^\infty(\Bbb R,\g) @()\L{T\evol_G }@(0,-1) @()\L{\cong}@(1,0) & 
C^\infty(\Bbb R,\g\rtimes\g) @()\l{\evol_{TG} }@(-1,-1) 
     @()\l{\evol_{\g\rtimes G}}@(0,-1) \\
TG @()\l{\cong}@(1,0) & \g\rtimes G
\endnewCD$$
For that we choose $X,Y\in C^\infty(\Bbb R,\g)$.
Let us first consider the evolution operator of the tangent group 
$TG$ in the picture $\g\rtimes G$. 
On $(\g,+)$ the evolution mapping is the definite integral, so 
going through the prescription 
\nmb!{5.8} for $\evol_{\g\rtimes G}$ we have in turn the 
following data:
$$\align
\evol&_{\g\rtimes G}(Y,X) = (h(1),g(1)),\quad\text{ where}\tag2\\
g(t):&=\Evol_G(X)(t) \in G,\\
Z(t):&= \Ad(g(t)\i).Y(t) \in \g,\\
h_0(t):&= \Evol_{(\g,+)}(Z)(t)
     =\int_0^t \Ad(g(u)\i).Y(u)\,du\in \g,\\
h(t):&= \Ad(g(t))h_0(t) = \Ad(g(t))\int_0^t \Ad(g(u)\i).Y(u)\,du\in \g.\\
\endalign$$
This shows the first equation in the corollary. 
The differential equation for the curve 
$(h(t),g(t))$, which by 
lemma \nmb!{5.2} has a unique solution starting at $(0,e)$, 
looks as follows, using \thetag1:
$$\align
\Bigl((h'(t),h(t)),g'(t)\Bigr) &= 
     T_{(0,e)}(\mu_{\g\rtimes G}^{(h(t),g(t))}).\Bigl((Y(t),0),X(t)\Bigr) \\
&= \Bigl(Y(t)+\Bigl(d\Ad(X(t)).h(t),0+\Ad(e).h(t)\Bigr),T(\mu_G^{g(t)}).X(t)\Bigr)\\
h'(t) &= Y(t) + \ad(X(t))h(t) \tag3\\
g'(t) &= T(\mu_G^{g(t)})X(t).
\endalign$$
For the computation of $T\evol_G$ we let
$$\gather 
g(t,s):= \evol_G\Bigl(u\mapsto t(X(tu)+sY(tu))\Bigr) = \Evol_G(X+sY)(t),\\
\text{ satisfying } \de^rg(\partial_t(t,s)) = X(t)+sY(t).
\endgather$$
Then $T\evol_G(Y,X)=\partial_s|_0g(1,s)$, 
and the derivative $\partial_s|_0 g(t,s)$ in $TG$ corresponds to
the element
$$
(T(\mu^{g(t,0)\i}).\partial_s|_0 g(t,s), g(t,0)) = 
(\de^rg(\partial_s(t,0)),g(t,0)) \in \g\rtimes G
$$
via right trivialization. For the right hand side we have 
$g(t,0)=g(t)$, so it remains to show that 
$\de^rg(\partial_s(t,0)) = h(t)$. We will show that 
$\de^rg(\partial_s(t,0))$ is the unique solution of the 
differential equation \thetag3 for $h(t)$.
Using the Maurer Cartan equation 
$d\de^rg - 
\tfrac12[\de^rg,\de^rg] = 0$
from lemma \nmb!{5.1} we get
$$\align
\partial_t\de^rg(\partial_s) &=
     \partial_s\de^rg(\partial_t) + 
     d(\de^rg)(\partial_t,\partial_s) + 
     \de^rg([\partial_t,\partial_s])\\
&= \partial_s\de^rg(\partial_t) + 
     [\de^rg(\partial_t),
     \de^rg(\partial_s)]_\g + 0\\
&= \partial_s(X(t)+sY(t)) + 
     [X(t)+sY(t),\de^rg(\partial_s)]_\g
\endalign$$
so that for $s=0$ we get 
$$\align
\partial_t\de^rg(\partial_s(t,0)) &=
     Y(t) + [X(t),\de^rg(\partial_s(t,0))]_\g \\
&=  Y(t) + \ad(X(t))\de^rg(\partial_s(t,0)).
\endalign$$
Thus $\de^rg(\partial_s(t,0))$ is a solution of the 
inhomogeneous linear ordinary differential equation \thetag3 as 
required.

It remains to check the last formula. Note that $X\mapsto tX(t\quad)$ 
is a bounded linear operator. So we have
$$\align
\Evol^r(X)(t) &= \evol(s\mapsto tX(ts)),\\
T_X(\Evol^r_G(\quad)&(t)).Y = T_{tX(t\quad)}\evol_G^r.(tY(t\quad))\\
&= T(\mu_{\evol_G^r(tX(t\quad))}).
     \int_0^1\Ad_G\Bigl(\Evol^r_G(tX(t\quad))(s)\i\Bigr).tY(ts)\,ds\\
&= T(\mu_{\Evol^r_G(X)(t)}).
     \int_0^1\Ad_G\Bigl(\evol^r_G(stX(st\quad))\i\Bigr).tY(ts)\,ds\\
&= T(\mu_{\Evol^r_G(X)(t)}).
     \int_0^t\Ad_G\Bigl(\Evol^r_G(X)(s)\i\Bigr).Y(s)\,ds.\qed
\endalign$$
\enddemo

\proclaim{\nmb.{5.10}. Corollary} 
For a regular Lie group $G$ the tangent mapping of the exponential mapping 
$\exp:\g\to G$ is given by:
$$\align
T_X\exp.Y &= T_e\mu_{\exp X}.\int_0^1 \Ad(\exp(-tX))Y dt\\
&= T_e\mu^{\exp X}.\int_0^1 \Ad(\exp(tX))Y dt
\endalign$$
\endproclaim

\remark{Remark}
This formula was first proved by {\rm\cit!{6}} for 
Lie groups with smooth exponential mapping. 
If $G$ is a Banach Lie group then we have from 
\nmb!{3.7}.\therosteritem4 and \nmb!{3.8} the series
$\Ad(\exp(tX))=\sum_{i=0}^\infty \tfrac{t^i}{i!}\ad(X)^i$, so that we get 
the usual formula
$$
T_X\exp = T_e\mu^{\exp X}.\sum_{i=0}^\infty \tfrac 1{(i+1)!}\ad(X)^i.
$$
\endremark 

\demo{Proof}
Just apply \nmb!{5.9} to constant curves $X,Y\in \g$.
\qed\enddemo

\head\totoc \nmb0{6}. Bundles with regular structure groups\endhead

\proclaim{\nmb.{6.1}. Theorem}
Let $(P,p,M,G)$ be a smooth (locally trivial) principal fiber bundle with a 
regular Lie group as structure group. Let $\om\in\Om^1(P,\g)$ be a 
principal connection form. 

Then the parallel
transport for the principal connection exists, is globally defined, 
and is $G$-equivariant. 
In detail: For each smooth curve $c:\Bbb R\to M$ there is a unique smooth
mapping $\Pt_c:\Bbb R\x P_{c(0)}\to P$ such that the following holds:
\roster
\item $\Pt(c,t,u)\in P_{c(t)}$, $\Pt(c,0) = Id_{P_{c(0)}}$, and
	$\om(\tfrac d{dt}\Pt(c,t,u))=0$.
\endroster
It has the following further properties:
\roster
\item [2]$\Pt(c,t):P_{c(0)}\to P_{c(t)}$ is $G$-equivariant, i\.e\.
	$\Pt(c,t,u.g) = \Pt(c,t,u).g$ holds for all $g\in G$ and $u\in P$. 
	Moreover we have $\Pt(c,t)^*(\ze_X|P_{c(t)})=\ze_X|P_{c(0)}$ for all 
	$X\in\frak g$.
\item For any smooth function $f:\Bbb R\to\Bbb R$ we have \newline
	$\Pt(c,f(t),u) = \Pt(c\o f,t,\Pt(c,f(0),u))$.
\item The parallel transport is smooth as a mapping
$$
\Pt: C^\infty(\Bbb R,M)\x_{(\ev_0,M,p\o \pr_2)}(\Bbb R\x P)\to P,
$$  
     where $C^\infty(\Bbb R,M)$ is considered as a smooth 
     space, see \cit!{4}, 1.4.1. 
\endroster
\endproclaim 

\demo{Proof}  
For a principal bundle chart $(U_\al,\ph_\al)$ we have the data from 
\nmb!{4.4}
$$\align
s_\al(x) :&= \ph_\al\i(x,e),\\
\om_\al :&= s_\al^*\om,\\
\om\o T(\ph_\al\i)&=(\ph_\al\i)^*\om\in\Om^1(U_\al\x G;\g)\\
(\ph_\al\i)^*\om(\xi_x,T\mu_g.X) &= (\ph_\al\i)^*\om(\xi_x,0_g) + X 
     =\Ad(g\i)\om_\al(\xi_x) + X.\\
\endalign$$
For a smooth curve $c:\Bbb R\to M$ the horizontal lift 
$\Pt(c,\quad,u)$ through $u\in P_{c(0)}$ is given by the ordinary 
differential equation $\om(\tfrac d{dt}\Pt(c,t,u))=0$ with initial 
condition $\Pt(c,0,u)=u$, among all smooth lifts of $c$. 
Locally we have
$$\ph_\al(\Pt(c,t,u))=(c(t),\ga(t)),$$
so that
$$\align
0&=\Ad(\ga(t))\om(\tfrac d{dt}\Pt(c,t,u)) = 
     \Ad(\ga(t))(\om\o T(\ph_\al\i))(c'(t),\ga'(t))\\
&=\Ad(\ga(t))((\ph_\al\i)^*\om)(c'(t),\ga'(t)) 
     = \om_\al(c'(t)) + T(\mu^{\ga(t)\i})\ga'(t),
\endalign$$
i\.e\. $\ga'(t) = - T(\mu^{\ga(t)}).\om_\al(c'(t))$, thus $\ga(t)$ is 
given by 
$$
\ga(t)=\Evol_G(-\om_\al(c'))(t).\ga(0)
     =\evol_G(s\mapsto -t\om_\al(c'(ts))).\ga(0).  
$$ 
By lemma \nmb!{5.2} we may glue the local solutions over different 
bundle charts $U_\al$, so $\Pt$ exists globally. 

Properties \therosteritem1 and \therosteritem3 are now clear,
and \therosteritem2 can be checked as follows: The condition 
$\om(\tfrac d{dt}\Pt(c,t,u).g)= \Ad(g\i)\om(\tfrac d{dt}\Pt(c,t,u))=0$ 
implies 
$\Pt(c,t,u).g=\Pt(c,t,u.g)$.
For the second assertion we compute for $u\in P_{c(0)}$:
$$\align
\Pt(c,t)^*(\ze_X|P_{c(t)})(u)&= T\Pt(c,t)\i\ze_X(\Pt(c,t,u)) \\
&= T\Pt(c,t)\i\tfrac d{ds}|_0\Pt(c,t,u).\exp(sX) \\
&= T\Pt(c,t)\i\tfrac d{ds}|_0\Pt(c,t,u.\exp(sX)) \\
&= \tfrac d{ds}|_0\Pt(c,t)\i\Pt(c,t,u.\exp(sX))\\
&= \tfrac d{ds}|_0 u.\exp(sX) = \ze_X(u).
\endalign$$

\therosteritem4. 
It suffices to check that $\Pt$ respects smooth curves. 
So let 
$(f,g):\Bbb R\to C^\infty(\Bbb R,M)\x_MP\subset 
     C^\infty(\Bbb R,M)\x P$
be a smooth curve.
By cartesian closedness of smooth spaces (see \cit!{4}, 1.4.3) 
the smooth curve $f:\Bbb R\to C^\infty(\Bbb R,M)$ 
corresponds to a smooth mapping 
$\hat f\in C^\infty(\Bbb R^2,M)$.
For a principal bundle chart $(U_\al,\ph_\al)$ as above we 
have $\ph_\al(\Pt(f(s),t,g(s)))=(f(s)(t),\ga(s,t))$, where $\ga$ is 
the evolution curve 
$$
\ga(s,t)=\Evol_G\Bigl(-\om_\al(\tfrac{\partial}{\partial t}
     \hat f(s,\quad))\Bigr)(t).\ph_\al(g(s)),
$$
which is clearly smooth in $(s,t)$.
\qed\enddemo

\proclaim{\nmb.{6.2}. Theorem}
Let $(P,p,M,G)$ be a smooth principal bundle with a 
regular Lie group as structure group. Let $\om\in\Om^1(P,\g)$ be a 
principal connection form. 
If the connection is flat, then the horizontal subbundle 
$H^\om(P):=\ker(\om)\subset TP$ is integrable and defines a foliation.

If $M$ is connected then each leaf of this horizontal foliation is a
covering of $M$. All leafs are isomorphic. 
\endproclaim

By standard arguments it follows that 
the principal bundle $P$ is associated to the 
universal covering of $M$ viewed as a principal fiber bundle with 
structure group the (discrete) fundamental group $\pi_1(M)$.

\demo{Proof}
Let $(U_\al,u_\al:U_\al\to u_\al(U_\al)\subset E_\al)$ 
be a smooth chart of the manifold 
$M$ and let $x_\al\in U_\al$ be such that $u_\al(x_\al)=0$ and the 
$c^\infty$-open subset $u_\al(U_\al)$ is disked in $E_\al$.
Let us also suppose that we have a principal fiber bundle chart 
$(U_\al,\ph_\al:P|U_\al\to U_\al\x G)$. 
We may cover $M$ by such $U_\al$.

We shall now construct for each $w_\al\in P_{x_\al}$ a smooth section
$\ps_\al:U_\al\to P$ whose image is an integral submanifold for the 
horizontal subbundle $\ker(\om)$. Namely, for $x\in U_\al$ let 
$c_x(t):=u_\al\i(tu_\al(x))$ for $t\in[0,1]$. Then we put 
$$
\ps_\al(x) := \Pt(c_x,1,w_\al).
$$
We have to show that the image of $T\ps_\al$ is contained in the 
horizontal bundle $\ker(\om)$. Then we get
$T_x\ps_\al = Tp|H^\om(p)_{\ps_\al(x)}\i$.
This is a consequence of the following notationally 
more suitable claim. 

Let $h:\Bbb R^2\to U_\al$ be smooth with $h(0,s)=x_\al$ for all $s$.
\newline
{\bf Claim:} $\dd s |Pt(h(.,s),1,w_\al)$ is horizontal.
 
Let $\ph_\al(w_\al)=(x_\al,g_\al)\in U_\al\x G$. 
Then from the proof of theorem \nmb!{6.1} we know that 
$$\align
\ph_\al\Pt(h(\quad,s),1,w_\al)&=(h(1,s),\ga(1,s)),\quad\text{ where }\\
\ga(t,s) &= \tilde\ga(t,s).g_\al\\
\tilde\ga(t,s)&=\evol_G\Bigl(u\mapsto-t\om_\al(\dd th(tu,s))\Bigr)\\
&= \Evol_G\Bigl(-(h^*\om_\al)(\partial_t(\quad,s))\Bigr)(t),\\
\om_\al &= s_\al^*\om,\qquad s_\al(x) = \ph_\al\i(x,e).
\endalign$$
Since the curvature $\Om=d\om+\tfrac12[\om,\om]_\wedge =0$ we have
$$\align
\partial_s(h^*\om_\al)(\partial_t) &= 
     \partial_t(h^*\om_\al)(\partial_s) - 
     d(h^*\om_\al)(\partial_t,\partial_s) - 
     (h^*\om_\al)([\partial_t,\partial_s])\\
&= \partial_t(h^*\om_\al)(\partial_s) + 
     [(h^*\om_\al)(\partial_t),(h^*\om_\al)(\partial_s)]_\g - 0.
\endalign$$
Using this and the expression for $T\evol_G$ from \nmb!{5.9} we 
have then:
$$\align
\dd s\tilde\ga(1,s) &= T_{-(h^*\om_\al)(\partial_t)(\quad,s)}\evol_G.
     \Bigl(-\partial_s(h^*\om_\al)(\partial_t)(\quad,s)\Bigr)\\
&= - T(\mu_{\tilde\ga(1,s)}).
     \int_0^1\Ad(\tilde\ga(t,s)\i)\partial_s(h^*\om_\al)(\partial_t)\,dt\\
&= - T(\mu_{\tilde\ga(1,s)}).\biggl(
     \int_0^1\Ad(\tilde\ga(t,s)\i)\partial_t(h^*\om_\al)(\partial_s)\,dt+\\
&\qquad\qquad  +\int_0^1\Ad(\tilde\ga(t,s)\i).\ad((h^*\om_\al)(\partial_t)).
     (h^*\om_\al)(\partial_s)\,dt\biggr).\\
\endalign$$
Next we integrate by parts, use \nmb!{3.9}.(3), and 
$\ka^l(\partial_t\tilde\ga(t,s)\i)=(h^*\om_\al)(\partial_t)(t,s)$ 
which follows from \nmb!{5.2}.
$$\align
\int_0^1&\Ad(\tilde\ga(t,s)\i)\partial_t(h^*\om_\al)(\partial_s)\,dt =\\
&= - \int_0^1\Bigl(\partial_t\Ad(\tilde\ga(t,s)\i)\Bigr)
     (h^*\om_\al)(\partial_s)\,dt
+ \Ad(\tilde\ga(t,s)\i)(h^*\om_\al)(\partial_s)\biggr|_{t=0}^{t=1}\\
&= - \int_0^1\Ad(\tilde\ga(t,s)\i).
     \ad\Bigl(\ka^l\partial_t(\tilde\ga(t,s)\i)\Bigr).
     (h^*\om_\al)(\partial_s)\,dt\\
&\quad + \Ad(\tilde\ga(1,s)\i)(h^*\om_\al)(\partial_s)(1,s) - 0\\
&= - \int_0^1\Ad(\tilde\ga(t,s)\i).\ad\Bigl((h^*\om_\al)(\partial_t)\Bigr).
     (h^*\om_\al)(\partial_s)\,dt\\
&\quad + \Ad(\tilde\ga(1,s)\i)(h^*\om_\al)(\partial_s)(1,s),
\endalign$$
so that finally 
$$\align
\dd s\tilde\ga(1,s) &= - T(\mu_{\tilde\ga(1,s)}).
     \Ad(\tilde\ga(1,s)\i)(h^*\om_\al)(\partial_s)(1,s)\\
&= - T(\mu^{\tilde\ga(1,s)}).(h^*\om_\al)(\partial_s)(1,s),\\
\dd s\ga(1,s) &= T(\mu^{g_\al}).\dd s\tilde\ga(1,s)\\
&= - T(\mu_{\ga(1,s)}).
     \Ad(\ga(1,s)\i)(h^*\om_\al)(\partial_s)(1,s)\\
\om(\dd s\Pt&(h(\quad,s),1,w_\al)) 
     = ((\ph_\al\i)^*\om)\Bigl(\dd sh(1,s),\dd s\ga(1,s)\Bigr)\\
&= \Ad(\ga(1,s)\i)\om_\al(\dd s h(1,s)) - 
\Ad(\ga(1,s)\i)(h^*\om_\al)(\partial_s)(1,s) =0,
\endalign$$
where in the end we used \nmb!{4.4}.(6). So the claim follows.

By the claim and by the uniqueness of parallel transport 
\nmb!{6.1}.\therosteritem1 for any smooth curve $c$ in $U_\al$ 
the horizontal curve $\ps_\al(c(t))$ 
coincides with $\Pt(c,t,\ps_\al(c(0)))$.
Moreover $U_\al\x G$ is $G$-equivariantly diffeomorphic to 
$p\i(U_\al)$ via $(x,g)\mapsto \ps_\al(x).g$. 

To finish the proof we may now glue overlapping right 
translations of $\ps_\al(U_\al)$ to maximal integral manifolds of 
the horizontal subbundle. As subset such an integral manifold 
consists of all endpoints of parallel transports of a fixed point. 
These are diffeomorphic covering spaces of $M$ via right 
translations.
\qed\enddemo

It is not clear, however, that the integral submanifolds of the 
theorem are initial submanifolds of $P$, or that they intersect each 
fiber in a totally disconnected subset, since $M$ might have 
uncountable fundamental group.

\subhead \nmb.{6.3}. Holonomy groups \endsubhead
Let $(P,p,M,G)$ be a principal fiber
bundle with regular structure group $G$ so that all parallel 
transports exist by theorem \nmb!{6.1}.
Let $\Ph = \ze\o\om$ be a principal connection. We assume
that $M$ is connected and we fix $x_0\in M$.

Now let us fix $u_0\in P_{x_0}$. 
Consider the subgroup $\operatorname{Hol}(\om,u_0)$ 
of the structure group $G$ which consists of 
all elements
$\tau(u_0,\Pt(c,t,u_0))\in G$ for $c$ any piecewise smooth closed
loop through $x_0$. Reparametrizing $c$ by a function which
is flat at each corner of $c$ we may assume that any $c$ is smooth.
We call $\operatorname{Hol}(\om,u_0)$ 
the \idx{\it holonomy group} of the connection.
If we consider only those curves $c$ which are nullhomotopic, we
obtain the \idx{\it restricted holonomy group} 
$\operatorname{Hol}_0(\om,x_0)$, 
a normal subgroup in $\operatorname{Hol}(\om,u_0)$.

\proclaim{ Theorem} 1. We have
$\operatorname{Hol}(\om,u_0.g)=\conj(g\i)\operatorname{Hol}(\om,u_0)$
and \newline
$\operatorname{Hol}_0(\om,u_0.g)=\conj(g\i)\operatorname{Hol}_0(\om,u_0)$.

2. For each curve $c$ in $M$ with $c(0)=x_0$ we have 
$\operatorname{Hol}(\om,\Pt(c,t,u_0))=\operatorname{Hol}(\om,u_0)$ 
and $\operatorname{Hol}_0(\om,\Pt(c,t,u_0))=\operatorname{Hol}_0(\om,u_0)$.
\endproclaim

\demo{Proof} 1. This follows from the properties of the mapping
$\tau$ from \nmb!{4.1} and from the from the $G$-equivariance of
the parallel transport: 
$$\ta(u_0.g,\Pt(c,1,u_0.g))= \ta(u_0,\Pt(c,1,u_0).g) 
=g\i.\ta(u_0,\Pt(c,1,u_0)).g.$$

2. By reparameterizing the curve $c$ we may assume that $t=1$, and we 
put $\Pt(c,1,u_0)=:u_1$. Then by definition for an element $g\in G$ 
we have $g\in\operatorname{Hol}(\om,u_1)$ if and only if 
$g=\ta(u_1,\Pt(e,1,u_1))$ for some closed smooth loop $e$ through 
$x_1:=c(1)=p(u_1)$, i\. e\. 
$$\align
\Pt(c,1)(r^g(u_0))&=r^g(\Pt(c,1)(u_0))=u_1g=\Pt(e,1)(\Pt(c,1)(u_0))\\
u_0g&=\Pt(c,1)\i\Pt(e,1)\Pt(c,1)(u_0)=\Pt(c.e.c\i,3)(u_0),
\endalign$$
where $c.e.c\i$ is the curve travelling along $c(t)$ for $0\le t\le 1$, 
along  $e(t-1)$ for $1\le t\le 3$, and along $c(3-t)$ for 
$2\le t\le 3$. This is equivalent to 
$g\in\operatorname{Hol}(\om,u_0)$.
Furthermore $e$ is nullhomotopic if and only if $c.e.c\i$ is 
nullhomotopic, so we also have 
$\operatorname{Hol}_0(\om,u_1)=\operatorname{Hol}_0(\om,u_0)$.
\qed\enddemo

\head\totoc\nmb0{7}. Rudiments of Lie theory for regular Lie groups 
\endhead

\subhead\nmb.{7.1}. From Lie algebras to Lie groups \endsubhead
It is not true in general that every convenient Lie 
agebra is the Lie  
algebra of a convenient Lie group. This is wrong for Banach Lie 
algebras and Banach Lie groups, one of the first examples is from 
\cit!{3}, see also \cit!{7}.

To Lie subalgebras in the Lie algebra of a Lie group do not 
correspond Lie subgroups in general, see the following easy example:

Let $\g\subset \X_c(\Bbb R^2)$ be the closed Lie subalgebra of all 
vector fields with compact support on $\Bbb R^2$ of the form 
$X(x,y)=f(x,y)\dd x +g(x,y)\dd y$ where $g$ vanishes on the strip 
$0\le x\le1$. 

\remark {Claim}
There is no Lie subgroup $G$ of $\Diff(\Bbb R^2)$ 
corresponding to $\g$.

If $G$ exists there is a smooth curve 
$t\mapsto f_t\in G\subset \Diff_c(\Bbb R^2)$ such that the smooth 
curve $X_t:=(\dd t 
f_t)\o f_t\i$ in $\g$ has the property that 
$X_0=f\dd x$ where $f=1$ near 0. But then $f_t$ moves the strip to 
the right for small $t$, so $\g$ is not invariant under
$\Ad^G(f_t)=f_t^*$, a contradiction. 

So we see that on any manifold of dimension greater that 2 there are 
closed Lie subalgebras of the Lie algebra of vector fields with 
compact support, which do not admit Lie subgroups.

Note that this example does not work for the Lie group of real 
analytic diffeomorphisms on a compact manifold, see \cit!{9}.
\endremark
 
\subhead\nmb.{7.2} \endsubhead
Let $G$ be a connected Lie group with 
Lie algebra $\g$.
For a smooth mapping $f:M\to G$ we considered in \nmb!{5.1} 
the right logarithmic derivative
$\de^r f\in\Om^1(M;\g)$ which 
is given by 
$\de^r f_x:=T(\mu^{f(x)\i})\o T_xf:T_xM\to T_{f(x)}G\to \g$
and which satisfies the 
left (from the left action) Maurer-Cartan equation
$$
d\de^rf - \frac12 
[\de^rf,\de^rf]^\g_\wedge=0.
$$
Similarly the 
left logarithmic derivative
$\de^l f\in\Om^1(M;\g)$ of $f\in C^\infty(M,G)$ was given by
$\de^l f_x := T(\mu_{f(x)\i})\o T_xf:T_xM\to T_{f(x)}G\to \g$
and satisfies the right (from the right action) Maurer Cartan 
equation 
$$
d\de^lf + \frac12 
[\de^lf,\de^lf]^\g_\wedge=0.
$$
For regular Lie groups we have the following converse:

\proclaim{Theorem} 
Let $G$ be a connected regular Lie group with 
Lie algebra $\g$.

If a 1-form $\ph\in\Om^1(M;\g)$ satisfies 
$d\ph-\frac12[\ph,\ph]_\wedge=0$ then for each simply connected subset 
$U\subset M$ there exists a smooth mapping $f:U\to G$ with 
$\de^r f=\ph|U$, and $f$ is uniquely detemined up to a right translation 
in $G$.

If a 1-form $\ps\in\Om^1(M;\g)$ satisfies 
$d\ps+\frac12[\ps,\ps]_\wedge=0$ then for each simply connected subset 
$U\subset M$ there exists a smooth mapping $f:U\to G$ with 
$\de^l f=\ps|U$, and $f$ is uniquely determined up to a left translation 
in $G$.
\endproclaim

The mapping $f$ is called the \idx{\it left developping} of $\ph$, or 
the \idx{\it right developping} of $\ps$, respectively. 

\demo{Proof} Let us treat the right logarithmic derivative 
since it leads to a principal connection for a bundle with right 
principal action. For the left logarithmic derivative the proof is 
similar, with the changes described in the second part of the proof 
of \nmb!{5.1}.

We put ourselves into the situation of the proof of \nmb!{5.1}.
If we are given a 1-form $\ph\in\Om^1(M;\g)$ with 
$d\ph - \tfrac12[\ph,\ph]_\wedge=0$ then we consider the 
1-form $\om^r\in\Om^1(M\x G;\g)$, given by the analogon of 
\nmb!{5.1}.\thetag1,
$$
\om^r = \ka^l - (\Ad\o\operatorname{Inv}).\ph\tag1
$$
Then $\om^r$ is a principal connection form on $M\x G$, since it 
reproduces the generators in $\g$ of the fundamental vector fields 
for the principal right action, i\.e\. the left invariant vector 
fields, and $\om^r$ is $G$-equivariant:
$$\align
((\mu^g)^*\om^r)_h &= \om^r_{hg}\o (Id\x T(\mu^g)) = 
     T(\mu_{g\i.h\i}).T(\mu^g) - \Ad(g\i.h\i).\ph\\
&= \Ad(g\i).\om^r_h.
\endalign$$
The computation in \nmb!{5.1}.\thetag3 for $\ph$ instead of 
$\de^r f$ shows 
that this connection is flat. 
Since the structure group $G$ is regular, by theorem \nmb!{6.2} 
the horizontal bundle is integrable, 
and $\operatorname{pr_1}:M\x G\to M$, restricted to each horizontal 
leaf, is a covering. Thus it may be inverted over each simply 
connected subset $U\subset M$, and the inverse $(Id,f):U\to M\x G$ is 
unique up to the choice of the branch of the covering, and the choice 
of the leaf, i\.e\. $f$ is unique up to a right translation by an 
element of $G$. The beginning of the proof of \nmb!{5.1} then shows that 
$\de^r f=\ph|U$. 
\qed\enddemo

\proclaim{\nmb.{7.3}. Theorem}
Let $G$ and $H$ be Lie groups with
Lie algebras $\frak g$ and $\frak h$, respectively. Let
$f:\frak g\to \frak h$ be a bounded homomorphism of Lie algebras.
If $H$ is regular and if $G$ is simply connected then there exists a 
unique homomorphism $F:G\to H$ of Lie groups with $T_eF=f$.
\endproclaim

This theorem is the main result in \cit!{16}, 
the proof there uses related methods. 

\demo{Proof}
We consider the 1-form 
$$
\ps\in \Om^1(G;\h),\quad \ps:= f\o \ka^r, \quad
\ps_g(\xi_g) = f(T(\mu^{g\i}).\xi_g),
$$
where $\ka^r$ is the right Maurer Cartan form from \nmb!{5.1}.
It satisfies the left Maurer Cartan equation 
$$\align
d\ps - \tfrac12[\ps,\ps]_{\wedge}^\h &=
     d(f\o \ka^r) - \tfrac12[f\o \ka^r,f\o \ka^r]_{\wedge}^\h \\
&= f\o (d\ka^r - \tfrac12[\ka^r,\ka^r]_{\wedge}^\g) = 0, \\
\endalign$$
by \nmb!{5.1}.\thetag{2'}. But then we can use theorem \nmb!{7.2} 
to conclude that there exists a unique smooth mapping $F:G\to H$ with 
$F(e)=e$ and whose right logarithmic derivative satisfies 
$\de^rF=\ps$. For $g\in G$ we have $(\mu^g)^*\ps=\ps$, thus also 
$$
\de^r(F\o \mu^g) = \de^rF\o T(\mu^g) = (\mu^g)^*\ps = \ps.
$$
By uniqueness in theorem \nmb!{7.2} again the mappings 
$F\o \mu^g, F:G\to H$ differ only by right translation in $H$ by 
$(F\o \mu^g)(e)=F(g)$, so that $F\o \mu^g = \mu^{F(g)}\o F$, or 
$F(g.g_1)=F(g).F(g_1)$. 
This also implies $F(g).F(g\i)=F(g.g\i)=F(e)=e$, so that $F$ is the 
unique homomorphism of Lie groups we looked for. 
\qed\enddemo

\proclaim{\nmb.{7.4}. Theorem} For a regular Lie group $G$ we have
$$\gather
\evol^r(X).\evol^r(Y) = 
     \evol^r\Bigl(t\mapsto X(t)+\Ad_G(\Evol^r(X)(t)).Y(t)\Bigr),\\
\evol^r(X)\i = \evol^r\Bigl(t\mapsto -\Ad_G(\Evol^r(X)(t)\i).X(t)\Bigr),
\endgather$$
so that $\evol^r:C^\infty(\Bbb R,\g)\to G$ is a surjective smooth 
homomorphism of Lie groups, where on $C^\infty(\Bbb R,\g)$ we 
consider the operations
$$\align
(X*Y)(t) &= X(t)+\Ad_G(\Evol^r(X)(t)).Y(t),\\
X\i(t) &= -\Ad_G(\Evol^r(X)(t)\i).X(t).
\endalign$$
With this operations and with 0 as unit element 
$(C^\infty(\Bbb R,\g),*)$ becomes again a regular Lie group.
Its Lie algebra is $C^\infty(\Bbb R,\g)$ with bracket 
$$\align
[X,Y]_{C^\infty(\Bbb R,\g)}(t) 
&= \Bigl[\int_0^tX(s)\,ds,Y(t)\Bigr]_\g + 
     \Bigl[X(t),\int_0^tY(s)\,ds\Bigr]_\g \\
&= \dd t \Bigl[\int_0^tX(s)\,ds,\int_0^tY(s)\,ds\Bigr]_\g.
\endalign$$
Its evolution operator is given by 
$$\align
\evol_{(C^\infty(\Bbb R,\g),*)}(X) :&= \Ad_G(\evol_G(Y^s)).
       \int_0^1 \Ad_G(\Evol_G(Y^s)(v)\i).X(v)(s)\,dv,\\
Y^s(t):&=\int_0^sX(t)(u)du.
\endalign$$
\endproclaim

\demo{Proof}
For $X,Y\in C^\infty(\Bbb R,\g)$ we compute
$$\align
\dd t &\Bigl(\Evol^r(X)(t).\Evol^r(Y)(t)\Bigr) =\\
&= T(\mu^{\Evol^r(Y)(t)}).T(\mu^{\Evol^r(X)(t)}).X(t) 
     + T(\mu_{\Evol^r(X)(t)}).T(\mu^{\Evol^r(Y)(t)}).Y(t)\\
&= T(\mu^{\Evol^r(X)(t).\Evol^r(Y)(t)}).(X(t)+\Ad_G(\Evol^r(X)(t))Y(t)),
\endalign$$
which implies also
$$
\Evol^r(X).\Evol^r(Y) = \Evol^r(X*Y),\quad \Evol(X)\i=\Evol(X\i).
$$
Thus $\Evol:C^\infty(\Bbb R,\g)\to C^\infty(\Bbb R,G)$ is a group 
isomorphism onto the subgroup $\{c\in C^\infty(\Bbb R,G):c(0)=e\}$ of 
$C^\infty(\Bbb R,G)$ 
with the pointwise product, 
which, however, is only a smooth space, see \cit!{4}, 1.4.1. 
Nevertheless it 
follows that the product on $C^\infty(\Bbb R,\g)$ is associative.
It is clear that these operations are smooth,
so that the convenient vector space 
$C^\infty(\Bbb R,\g)$ becomes a Lie group; and $C^\infty(\Bbb R,G)$ 
becomes a manifold.

Now we aim for the Lie bracket.
We have
$$\align
(X*Y*X\i&)(t) = \biggl(\Bigl(X+\Ad(\Evol^r(X)).Y\Bigr)*
     \Bigl(-\Ad(\Evol^r(X)\i).X\Bigr)\biggr)(t)\\
&= X(t)+\Ad(\Evol^r(X)(t)).Y(t) -\\
&\quad  -\Ad\Bigl(\Evol^r(X*Y)(t)\Bigr).\Ad\Bigl(\Evol^r(X)(t)\i\Bigr).X(t)\\
&= X(t)+\Ad\Bigl(\Evol^r(X)(t)\Bigr).Y(t) -\\
&\quad  -\Ad\Bigl(\Evol^r(X)(t)\Bigr).
     \Ad\Bigl(\Evol^r(Y)(t)\Bigr).
     \Ad\Bigl(\Evol^r(X)(t)\i\Bigr).X(t).\\
\endalign$$
We shall need 
$$\align
T_0\Bigl(\Ad_G(\Evol^r(\quad)(t))\Bigr).Y &= 
     T_e\Ad_G.T_0(\Evol^r(\quad)(t)).Y\\
&= \ad_\g\Bigl(\int_0^tY(s)\,ds\Bigr),\quad\text{ by }\nmb!{5.9}.
\endalign$$
Using this we can differentiate the conjugation,
$$\align
(\Ad&_{C^\infty(\Bbb R,\g)}(X).Y)(t) = (T_0(X*(\quad)*X\i).Y)(t) \\
&= 0 + \Ad(\Evol^r(X)(t)).Y(t) -\\
&\quad  -\Ad(\Evol^r(X)(t)).\Bigl(T_0(\Ad(\Evol^r(\quad)(t))).Y\Bigr).
     \Ad(\Evol^r(X)(t)\i).X(t)\\
&= \Ad(\Evol^r(X)(t)).Y(t) -\\
&\quad  -\Ad(\Evol^r(X)(t)).\ad_\g\Bigl(\int_0^tY(s)\,ds\Bigr).
     \Ad(\Evol^r(X)(t)\i).X(t)\\
&= \Ad(\Evol^r(X)(t)).Y(t) -
     \ad_\g.\Bigl(\Ad(\Evol^r(X)(t)).\int_0^tY(s)\,ds\Bigr).X(t).\\
\endalign$$
Now we can compute the Lie bracket
$$\align
[X,Y&]_{C^\infty(\Bbb R,\g)}(t) 
     = \Bigl(T_0(\Ad_{C^\infty(\Bbb R,\g)}(\quad).Y).X\Bigr)(t)\\
&= T_0\Bigl(\Ad(\Evol^r(\quad)(t)).X\Bigr).Y(t) 
     -0 -\Bigl[\Ad(\Evol^r(0)(t)).\int_0^tY(s)\,ds,X(t)\Bigr]_\g\\
&= \Bigl[\int_0^tX(s)\,ds,Y(t)\Bigr]_\g - 
     \Bigl[\int_0^tY(s)\,ds,X(t)\Bigr]_\g\\
&= \Bigl[\int_0^tX(s)\,ds,Y(t)\Bigr]_\g + 
     \Bigl[X(t),\int_0^tY(s)\,ds\Bigr]_\g\\
&= \dd t \Bigl[\int_0^tX(s)\,ds,\int_0^tY(s)\,ds\Bigr]_\g.
\endalign$$
Now we show that the Lie group $(C^\infty(\Bbb R,\g),*)$ is regular.
Let $\check X\in C^\infty(\Bbb R,C^\infty(\Bbb R,\g))$ correspond to 
$X\in C^\infty(\Bbb R^2,\g)$. We look for 
$g\in C^\infty(\Bbb R^2,\g)$ which satisfies the equation 
\nmb!{5.3}.\thetag1:
$$\align
\mu^{g(t,\quad)}(Y)(s) &= (Y * g(t,\quad))(s) 
     = Y(s) + \Ad_G(\Evol_G(Y)(s)).g(t,s)\\
\dd t g(t,s) &= \Bigl(T_0(\mu^{g(t,\quad)}).X(t,\quad)\Bigr)(s)\\
&= X(t,s) + \Bigl(T_0\Bigl(\Ad_G(\Evol_G(\quad)(s))\Bigr).
     X(t,\quad)\Bigr).g(t,s)\\
&= X(t,s) + \ad_\g\Bigl(\int_0^sX(t,u)du\Bigr).g(t,s)\\
&= X(t,s) + \Bigl[\int_0^sX(t,u)du),g(t,s)\Bigr]_\g.
\endalign$$
This is the differential equation \nmb!{5.9}.\thetag3, depending 
smoothly on a further parameter $s$, which has the following 
unique solution which is given by \nmb!{5.9}.\thetag2
$$\align
g(t,s) :&= \Ad_G(\Evol_G(Y^s)(t)).
       \int_0^t \Ad_G(\Evol_G(Y^s)(v)\i).X(v,s)\,dv\\
Y^s(t):&=\int_0^sX(t,u)du.
\endalign$$
Since this solution is visibly smooth in $X$, the Lie group 
$C^\infty(\Bbb R,\g)$ is regular.
For convenience (yours, not ours) we show now (once more) that this 
is a solution. Putting $Y^s(t):=\int_0^sX(t,u)du$ we have by 
\nmb!{3.9}.\therosteritem3
$$\align
\dd t g&(t,s) =\\
&= d\Ad (\dd t\Evol(Y^s)(t)).
     \int_0^t \Ad(\Evol(Y^s)(v)\i).X(v,s)\,dv\\
&\quad+ \Ad(\Evol(Y^s)(t)). 
     \Ad(\Evol(Y^s)(t)\i).X(t,s)\\
\allowdisplaybreak
&= ((\ad\o \ka^r).\Ad)\Bigl(T(\mu^{\Evol(Y^s)(t)}).Y^s(t)\Bigr).
     \int_0^t \Ad(\Evol(Y^s)(v)\i).X(v,s)\,dv\\
&\quad + X(t,s)\\
&= \ad(Y^s(t)).\Ad(\Evol(Y^s)(t)).
     \int_0^t \Ad(\Evol(Y^s)(v)\i).X(v,s)\,dv + X(t,s)\\
&= \Bigl[\int_0^sX(t,u)du),g(t,s)\Bigr]_\g + X(t,s).\qed
\endalign$$
\enddemo

\proclaim{\nmb.{7.5}. Corollary} Let $G$ be a regular Lie group.
Then as smooth spaces and groups we have the following isomorphims 
$$
(C^\infty(\Bbb R,\g),*) \rtimes G 
\cong \{f\in C^\infty(\Bbb R,G):f(0)=e\} \rtimes G 
\cong C^\infty(\Bbb R,G),
$$
where $g\in G$ acts on $f$ by $(\al_g(f))(t)=g.f(t).g\i$, and on 
$X\in C^\infty(\Bbb R,\g)$ by $\al_g(X)(t)=\Ad_G(g)(X(t))$.
The leftmost space is a smooth manifold, thus all spaces are regular 
Lie groups.

For the Lie algebras we have an isomorphism
$$\gather
C^\infty(\Bbb R,\g)\rtimes \g \cong C^\infty(\Bbb R,\g),\\
(X,\et)\mapsto \Bigl(t\mapsto \et + \int_0^t X(s)ds\Bigr)\\
(Y',Y(0)) \gets Y
\endgather$$
where on the left hand side the Lie bracket is given by 
$$\multline
[(X_1,\et_1),(X_2,\et_2)] =\\
= \Bigl(t\mapsto [\tsize\int_0^tX_1(s)\,ds,X_2(t)]_\g + 
     [X_1(t),\tsize\int_0^tX_2(s)\,ds]_\g + [\et_1,X_2(t)]_\g - 
     [\et_2,X_1]_\g,\\ [\et_1,\et_2]_\g\Bigr),
\endmultline$$
and where on the right hand side the bracket is given by 
$$
[X,Y](t) = [X(t),Y(t)]_\g.
$$
On the right hand sides the evolution operator is given by 
$$\Evol^r_{C^\infty(\Bbb R,G)}= C^\infty(\Bbb R,\Evol^r_G).$$
\endproclaim

\subhead\nmb.{7.6}. Remarks  \endsubhead
Let $G$ be a connected regular Lie group.
The smooth homomorphism $\evol_G^r:C^\infty(\Bbb R,\g)\to G$ admits 
local smooth sections. Namely using a smooth chart near $e$ of $G$ we 
can choose a smooth curve $c_g:\Bbb R\to G$ with $c_g(0)=e$ and 
$c_g(1)=g$, depending smoothly on $g$, for $g$ near $e$. Then 
$s(g):=\de^rc_g$ is a local smooth section. We have an extension 
of groups
$$
0 \to K \to C^\infty(\Bbb R,\g) \East{\evol_G^r}{} G \to \{e\}
$$
where $K=\ker(\evol_G^r)$ is isomorphic 
to the smooth group $\{f\in C^\infty(\Bbb R,G):f(0)=e, f(1)=e\}$ via 
the mapping $\Evol_G^r$.
We do not know whether $K$ is a submanifold.

Next we consider the smooth 
group $C^\infty((S^1,1),(G,e))$ of all smooth mappings $f:S^1\to G$ 
with $f(1)=e$. With pointwise multiplication
this is a splitting closed normal subgroup of the regular Lie 
group $C^\infty(S^1,G)$ with the manifold structure described in 
\cit!{10} and \cit!{12}. 
Moreover $C^\infty(S^1,G)$ is the semidirect product 
$C^\infty((S^1,1),(G,e))\rtimes G$, where $G$ acts by conjugation on 
$C^\infty((S^1,1),(G,e))$. So by theorem \nmb!{5.5} the subgroup 
$C^\infty((S^1,1),(G,e))$ is also regular.

The right logarithmic derivative 
$\de^r:C^\infty(S^1,G)\to C^\infty(S^1,\g)$ restricts to a 
diffeomorphism 
$C^\infty((S^1,1),(G,e))\to \ker(\evol_G)\subset C^\infty(S^1,\g)$,
thus $\ker(\evol_G:C^\infty(S^1,\g)\to G)$ is a regular Lie group 
isomorphic to $C^\infty((S^1,1),(G,e))$. It is also a subgroup (via 
pullback by the covering mapping $e^{2\pi i t}:\Bbb R\to S^1$) of the 
regular Lie group $(C^\infty(\Bbb R,\g),*)$.  
Note that $C^\infty(S^1,\g)$ is not a subgroup, it is not closed 
under the product $*$, if $G$ is not abelian.

\enddocument

%% file: diag.tex
%
\catcode`\@=11

\def\East#1#2{\setboxz@h{$\m@th\ssize\;{#1}\;\;$}%
 \setbox@ne\hbox{$\m@th\ssize\;{#2}\;\;$}\setbox\tw@\hbox{$\m@th#2$}%
 \dimen@\minaw@
 \ifdim\wdz@>\dimen@ \dimen@\wdz@ \fi  \ifdim\wd@ne>\dimen@ \dimen@\wd@ne \fi
 \ifdim\wd\tw@>\z@
  \mathrel{\mathop{\hbox to\dimen@{\rightarrowfill}}\limits^{#1}_{#2}}%
 \else
  \mathrel{\mathop{\hbox to\dimen@{\rightarrowfill}}\limits^{#1}}%
 \fi}
\def\West#1#2{\setboxz@h{$\m@th\ssize\;\;{#1}\;$}%
 \setbox@ne\hbox{$\m@th\ssize\;\;{#2}\;$}\setbox\tw@\hbox{$\m@th#2$}%
 \dimen@\minaw@
 \ifdim\wdz@>\dimen@ \dimen@\wdz@ \fi \ifdim\wd@ne>\dimen@ \dimen@\wd@ne \fi
 \ifdim\wd\tw@>\z@
  \mathrel{\mathop{\hbox to\dimen@{\leftarrowfill}}\limits^{#1}_{#2}}%
 \else
  \mathrel{\mathop{\hbox to\dimen@{\leftarrowfill}}\limits^{#1}}%
 \fi}
\font\arrow@i=lams1
\font\arrow@ii=lams2
\font\arrow@iii=lams3
\font\arrow@iv=lams4
\font\arrow@v=lams5
\newbox\zer@
\newdimen\standardcgap
\standardcgap=40\p@
\newdimen\hunit
\hunit=\tw@\p@
\newdimen\standardrgap
\standardrgap=32\p@
\newdimen\vunit
\vunit=1.6\p@
\def\Cgaps#1{\RIfM@
  \standardcgap=#1\standardcgap\relax \hunit=#1\hunit\relax
 \else \nonmatherr@\Cgaps \fi}
\def\Rgaps#1{\RIfM@
  \standardrgap=#1\standardrgap\relax \vunit=#1\vunit\relax
 \else \nonmatherr@\Rgaps \fi}
\newdimen\getdim@
\def\getcgap@#1{\ifcase#1\or\getdim@\z@\else\getdim@\standardcgap\fi}
\def\getrgap@#1{\ifcase#1\getdim@\z@\else\getdim@\standardrgap\fi}
\def\cgaps#1{\RIfM@
 \cgaps@{#1}\edef\getcgap@##1{\i@=##1\relax\the\toks@}\toks@{}\else
 \nonmatherr@\cgaps\fi}
\def\rgaps#1{\RIfM@
 \rgaps@{#1}\edef\getrgap@##1{\i@=##1\relax\the\toks@}\toks@{}\else
 \nonmatherr@\rgaps\fi}
\def\Gaps@@{\gaps@@}
\def\cgaps@#1{\toks@{\ifcase\i@\or\getdim@=\z@}%
 \gaps@@\standardcgap#1;\gaps@@\gaps@@
 \edef\next@{\the\toks@\noexpand\else\noexpand\getdim@\noexpand\standardcgap
  \noexpand\fi}%
 \toks@=\expandafter{\next@}}
\def\rgaps@#1{\toks@{\ifcase\i@\getdim@=\z@}%
 \gaps@@\standardrgap#1;\gaps@@\gaps@@
 \edef\next@{\the\toks@\noexpand\else\noexpand\getdim@\noexpand\standardrgap
  \noexpand\fi}%
 \toks@=\expandafter{\next@}}
\def\gaps@@#1#2;#3{\mgaps@#1#2\mgaps@
 \edef\next@{\the\toks@\noexpand\or\noexpand\getdim@
  \noexpand#1\the\mgapstoks@@}%
 \global\toks@=\expandafter{\next@}%
 \DN@{#3}%
 \ifx\next@\Gaps@@\gdef\next@##1\gaps@@{}\else
  \gdef\next@{\gaps@@#1#3}\fi\next@}
\def\mgaps@#1{\let\mgapsnext@#1\FN@\mgaps@@}
\def\mgaps@@{\ifx\next\space@\DN@. {\FN@\mgaps@@}\else
 \DN@.{\FN@\mgaps@@@}\fi\next@.}
\def\mgaps@@@{\ifx\next\w\let\next@\mgaps@@@@\else
 \let\next@\mgaps@@@@@\fi\next@}
\newtoks\mgapstoks@@
\def\mgaps@@@@@#1\mgaps@{\getdim@\mgapsnext@\getdim@#1\getdim@
 \edef\next@{\noexpand\getdim@\the\getdim@}%
 \mgapstoks@@=\expandafter{\next@}}
\def\mgaps@@@@\w#1#2\mgaps@{\mgaps@@@@@#2\mgaps@
 \setbox\zer@\hbox{$\m@th\hskip15\p@\tsize@#1$}%
 \dimen@\wd\zer@
 \ifdim\dimen@>\getdim@ \getdim@\dimen@ \fi
 \edef\next@{\noexpand\getdim@\the\getdim@}%
 \mgapstoks@@=\expandafter{\next@}}
\def\changewidth#1#2{\setbox\zer@\hbox{$\m@th#2}%
 \hbox to\wd\zer@{\hss$\m@th#1$\hss}}
\atdef@({\FN@\ARROW@}
\def\ARROW@{\ifx\next)\let\next@\OPTIONS@\else
 \DN@{\csname\string @(\endcsname}\fi\next@}
\newif\ifoptions@
\def\OPTIONS@){\ifoptions@\let\next@\relax\else
 \DN@{\options@true\begingroup\optioncodes@}\fi\next@}
\newif\ifN@
\newif\ifE@
\newif\ifNESW@
\newif\ifH@
\newif\ifV@
\newif\ifHshort@
\expandafter\def\csname\string @(\endcsname #1,#2){%
 \ifoptions@\let\next@\endgroup\else\let\next@\relax\fi\next@
 \N@false\E@false\H@false\V@false\Hshort@false
 \ifnum#1>\z@\E@true\fi
 \ifnum#1=\z@\V@true\tX@false\tY@false\a@false\fi
 \ifnum#2>\z@\N@true\fi
 \ifnum#2=\z@\H@true\tX@false\tY@false\a@false\ifshort@\Hshort@true\fi\fi
 \NESW@false
 \ifN@\ifE@\NESW@true\fi\else\ifE@\else\NESW@true\fi\fi
 \arrow@{#1}{#2}%
 \global\options@false
 \global\scount@\z@\global\tcount@\z@\global\arrcount@\z@
 \global\s@false\global\sxdimen@\z@\global\sydimen@\z@
 \global\tX@false\global\tXdimen@i\z@\global\tXdimen@ii\z@
 \global\tY@false\global\tYdimen@i\z@\global\tYdimen@ii\z@
 \global\a@false\global\exacount@\z@
 \global\x@false\global\xdimen@\z@
 \global\X@false\global\Xdimen@\z@
 \global\y@false\global\ydimen@\z@
 \global\Y@false\global\Ydimen@\z@
 \global\p@false\global\pdimen@\z@
 \global\label@ifalse\global\label@iifalse
 \global\dl@ifalse\global\ldimen@i\z@
 \global\dl@iifalse\global\ldimen@ii\z@
 \global\short@false\global\unshort@false}
\newif\iflabel@i
\newif\iflabel@ii
\newcount\scount@
\newcount\tcount@
\newcount\arrcount@
\newif\ifs@
\newdimen\sxdimen@
\newdimen\sydimen@
\newif\iftX@
\newdimen\tXdimen@i
\newdimen\tXdimen@ii
\newif\iftY@
\newdimen\tYdimen@i
\newdimen\tYdimen@ii
\newif\ifa@
\newcount\exacount@
\newif\ifx@
\newdimen\xdimen@
\newif\ifX@
\newdimen\Xdimen@
\newif\ify@
\newdimen\ydimen@
\newif\ifY@
\newdimen\Ydimen@
\newif\ifp@
\newdimen\pdimen@
\newif\ifdl@i
\newif\ifdl@ii
\newdimen\ldimen@i
\newdimen\ldimen@ii
\newif\ifshort@
\newif\ifunshort@
\def\zero@#1{\ifnum\scount@=\z@
 \if#1e\global\scount@\m@ne\else
 \if#1t\global\scount@\tw@\else
 \if#1h\global\scount@\thr@@\else
 \if#1'\global\scount@6 \else
 \if#1`\global\scount@7 \else
 \if#1(\global\scount@8 \else
 \if#1)\global\scount@9 \else
 \if#1s\global\scount@12 \else
 \if#1H\global\scount@13 \else
 \Err@{\Invalid@@ option \string\0}\fi\fi\fi\fi\fi\fi\fi\fi\fi
 \fi}
\def\one@#1{\ifnum\tcount@=\z@
 \if#1e\global\tcount@\m@ne\else
 \if#1h\global\tcount@\tw@\else
 \if#1t\global\tcount@\thr@@\else
 \if#1'\global\tcount@4 \else
 \if#1`\global\tcount@5 \else
 \if#1(\global\tcount@10 \else
 \if#1)\global\tcount@11 \else
 \if#1s\global\tcount@12 \else
 \if#1H\global\tcount@13 \else
 \Err@{\Invalid@@ option \string\1}\fi\fi\fi\fi\fi\fi\fi\fi\fi
 \fi}
\def\a@#1{\ifnum\arrcount@=\z@
 \if#10\global\arrcount@\m@ne\else
 \if#1+\global\arrcount@\@ne\else
 \if#1-\global\arrcount@\tw@\else
 \if#1=\global\arrcount@\thr@@\else
 \Err@{\Invalid@@ option \string\a}\fi\fi\fi\fi
 \fi}
\def\ds@(#1;#2){\ifs@\else
 \global\s@true
 \sxdimen@\hunit \global\sxdimen@#1\sxdimen@\relax
 \sydimen@\vunit \global\sydimen@#2\sydimen@\relax
 \fi}
\def\dtX@(#1;#2){\iftX@\else
 \global\tX@true
 \tXdimen@i\hunit \global\tXdimen@i#1\tXdimen@i\relax
 \tXdimen@ii\vunit \global\tXdimen@ii#2\tXdimen@ii\relax
 \fi}
\def\dtY@(#1;#2){\iftY@\else
 \global\tY@true
 \tYdimen@i\hunit \global\tYdimen@i#1\tYdimen@i\relax
 \tYdimen@ii\vunit \global\tYdimen@ii#2\tYdimen@ii\relax
 \fi}
\def\da@#1{\ifa@\else\global\a@true\global\exacount@#1\relax\fi}
\def\dx@#1{\ifx@\else
 \global\x@true
 \xdimen@\hunit \global\xdimen@#1\xdimen@\relax
 \fi}
\def\dX@#1{\ifX@\else
 \global\X@true
 \Xdimen@\hunit \global\Xdimen@#1\Xdimen@\relax
 \fi}
\def\dy@#1{\ify@\else
 \global\y@true
 \ydimen@\vunit \global\ydimen@#1\ydimen@\relax
 \fi}
\def\dY@#1{\ifY@\else
 \global\Y@true
 \Ydimen@\vunit \global\Ydimen@#1\Ydimen@\relax
 \fi}
\def\p@@#1{\ifp@\else
 \global\p@true
 \pdimen@\hunit \divide\pdimen@\tw@ \global\pdimen@#1\pdimen@\relax
 \fi}
\def\L@#1{\iflabel@i\else
 \global\label@itrue \gdef\label@i{#1}%
 \fi}
\def\l@#1{\iflabel@ii\else
 \global\label@iitrue \gdef\label@ii{#1}%
 \fi}
\def\dL@#1{\ifdl@i\else
 \global\dl@itrue \ldimen@i\hunit \global\ldimen@i#1\ldimen@i\relax
 \fi}
\def\dl@#1{\ifdl@ii\else
 \global\dl@iitrue \ldimen@ii\hunit \global\ldimen@ii#1\ldimen@ii\relax
 \fi}
\def\s@{\ifunshort@\else\global\short@true\fi}
\def\uns@{\ifshort@\else\global\unshort@true\global\short@false\fi}
\def\optioncodes@{\let\0\zero@\let\1\one@\let\a\a@\let\ds\ds@\let\dtX\dtX@
 \let\dtY\dtY@\let\da\da@\let\dx\dx@\let\dX\dX@\let\dY\dY@\let\dy\dy@
 \let\p\p@@\let\L\L@\let\l\l@\let\dL\dL@\let\dl\dl@\let\s\s@\let\uns\uns@}
\def\slopes@{\\161\\152\\143\\134\\255\\126\\357\\238\\349\\45{10}\\56{11}%
 \\11{12}\\65{13}\\54{14}\\43{15}\\32{16}\\53{17}\\21{18}\\52{19}\\31{20}%
 \\41{21}\\51{22}\\61{23}}
\newcount\tan@i
\newcount\tan@ip
\newcount\tan@ii
\newcount\tan@iip
\newdimen\slope@i
\newdimen\slope@ip
\newdimen\slope@ii
\newdimen\slope@iip
\newcount\angcount@
\newcount\extracount@
\def\slope@{{\slope@i=\secondy@ \advance\slope@i-\firsty@
 \ifN@\else\multiply\slope@i\m@ne\fi
 \slope@ii=\secondx@ \advance\slope@ii-\firstx@
 \ifE@\else\multiply\slope@ii\m@ne\fi
 \ifdim\slope@ii<\z@
  \global\tan@i6 \global\tan@ii\@ne \global\angcount@23
 \else
  \dimen@\slope@i \multiply\dimen@6
  \ifdim\dimen@<\slope@ii
   \global\tan@i\@ne \global\tan@ii6 \global\angcount@\@ne
  \else
   \dimen@\slope@ii \multiply\dimen@6
   \ifdim\dimen@<\slope@i
    \global\tan@i6 \global\tan@ii\@ne \global\angcount@23
   \else
    \tan@ip\z@ \tan@iip \@ne
    \def\\##1##2##3{\global\angcount@=##3\relax
     \slope@ip\slope@i \slope@iip\slope@ii
     \multiply\slope@iip##1\relax \multiply\slope@ip##2\relax
     \ifdim\slope@iip<\slope@ip
      \global\tan@ip=##1\relax \global\tan@iip=##2\relax
     \else
      \global\tan@i=##1\relax \global\tan@ii=##2\relax
      \def\\####1####2####3{}%
     \fi}%
    \slopes@
    \slope@i=\secondy@ \advance\slope@i-\firsty@
    \ifN@\else\multiply\slope@i\m@ne\fi
    \multiply\slope@i\tan@ii \multiply\slope@i\tan@iip \multiply\slope@i\tw@
    \count@\tan@i \multiply\count@\tan@iip
    \extracount@\tan@ip \multiply\extracount@\tan@ii
    \advance\count@\extracount@
    \slope@ii=\secondx@ \advance\slope@ii-\firstx@
    \ifE@\else\multiply\slope@ii\m@ne\fi
    \multiply\slope@ii\count@
    \ifdim\slope@i<\slope@ii
     \global\tan@i=\tan@ip \global\tan@ii=\tan@iip
     \global\advance\angcount@\m@ne
    \fi
   \fi
  \fi
 \fi}%
}
\def\slope@a#1{{\def\\##1##2##3{\ifnum##3=#1\global\tan@i=##1\relax
 \global\tan@ii=##2\relax\fi}\slopes@}}
\newcount\i@
\newcount\j@
\newcount\colcount@
\newcount\Colcount@
\newcount\tcolcount@
\newdimen\rowht@
\newdimen\rowdp@
\newcount\rowcount@
\newcount\Rowcount@
\newcount\maxcolrow@
\newtoks\colwidthtoks@
\newtoks\Rowheighttoks@
\newtoks\Rowdepthtoks@
\newtoks\widthtoks@
\newtoks\Widthtoks@
\newtoks\heighttoks@
\newtoks\Heighttoks@
\newtoks\depthtoks@
\newtoks\Depthtoks@
\newif\iffirstnewCDcr@
\def\dotoks@i{%
 \global\widthtoks@=\expandafter{\the\widthtoks@\else\getdim@\z@\fi}%
 \global\heighttoks@=\expandafter{\the\heighttoks@\else\getdim@\z@\fi}%
 \global\depthtoks@=\expandafter{\the\depthtoks@\else\getdim@\z@\fi}}
\def\dotoks@ii{%
 \global\widthtoks@{\ifcase\j@}%
 \global\heighttoks@{\ifcase\j@}%
 \global\depthtoks@{\ifcase\j@}}
\def\prenewCD@#1\endnewCD{\setbox\zer@
 \vbox{%
  \def\arrow@##1##2{{}}%
  \rowcount@\m@ne \colcount@\z@ \Colcount@\z@
  \firstnewCDcr@true \toks@{}%
  \widthtoks@{\ifcase\j@}%
  \Widthtoks@{\ifcase\i@}%
  \heighttoks@{\ifcase\j@}%
  \Heighttoks@{\ifcase\i@}%
  \depthtoks@{\ifcase\j@}%
  \Depthtoks@{\ifcase\i@}%
  \Rowheighttoks@{\ifcase\i@}%
  \Rowdepthtoks@{\ifcase\i@}%
  \Let@
  \everycr{%
   \noalign{%
    \global\advance\rowcount@\@ne
    \ifnum\colcount@<\Colcount@
    \else
     \global\Colcount@=\colcount@ \global\maxcolrow@=\rowcount@
    \fi
    \global\colcount@\z@
    \iffirstnewCDcr@
     \global\firstnewCDcr@false
    \else
     \edef\next@{\the\Rowheighttoks@\noexpand\or\noexpand\getdim@\the\rowht@}%
      \global\Rowheighttoks@=\expandafter{\next@}%
     \edef\next@{\the\Rowdepthtoks@\noexpand\or\noexpand\getdim@\the\rowdp@}%
      \global\Rowdepthtoks@=\expandafter{\next@}%
     \global\rowht@\z@ \global\rowdp@\z@
     \dotoks@i
     \edef\next@{\the\Widthtoks@\noexpand\or\the\widthtoks@}%
      \global\Widthtoks@=\expandafter{\next@}%
     \edef\next@{\the\Heighttoks@\noexpand\or\the\heighttoks@}%
      \global\Heighttoks@=\expandafter{\next@}%
     \edef\next@{\the\Depthtoks@\noexpand\or\the\depthtoks@}%
      \global\Depthtoks@=\expandafter{\next@}%
     \dotoks@ii
    \fi}}%
  \tabskip\z@
  \halign{&\setbox\zer@\hbox{\vrule height10\p@ width\z@ depth\z@
   $\m@th\displaystyle{##}$}\copy\zer@
   \ifdim\ht\zer@>\rowht@ \global\rowht@\ht\zer@ \fi
   \ifdim\dp\zer@>\rowdp@ \global\rowdp@\dp\zer@ \fi
   \global\advance\colcount@\@ne
   \edef\next@{\the\widthtoks@\noexpand\or\noexpand\getdim@\the\wd\zer@}%
    \global\widthtoks@=\expandafter{\next@}%
   \edef\next@{\the\heighttoks@\noexpand\or\noexpand\getdim@\the\ht\zer@}%
    \global\heighttoks@=\expandafter{\next@}%
   \edef\next@{\the\depthtoks@\noexpand\or\noexpand\getdim@\the\dp\zer@}%
    \global\depthtoks@=\expandafter{\next@}%
   \cr#1\crcr}}%
 \Rowcount@=\rowcount@
 \global\Widthtoks@=\expandafter{\the\Widthtoks@\fi\relax}%
 \edef\Width@##1##2{\i@=##1\relax\j@=##2\relax\the\Widthtoks@}%
 \global\Heighttoks@=\expandafter{\the\Heighttoks@\fi\relax}%
 \edef\Height@##1##2{\i@=##1\relax\j@=##2\relax\the\Heighttoks@}%
 \global\Depthtoks@=\expandafter{\the\Depthtoks@\fi\relax}%
 \edef\Depth@##1##2{\i@=##1\relax\j@=##2\relax\the\Depthtoks@}%
 \edef\next@{\the\Rowheighttoks@\noexpand\fi\relax}%
 \global\Rowheighttoks@=\expandafter{\next@}%
 \edef\Rowheight@##1{\i@=##1\relax\the\Rowheighttoks@}%
 \edef\next@{\the\Rowdepthtoks@\noexpand\fi\relax}%
 \global\Rowdepthtoks@=\expandafter{\next@}%
 \edef\Rowdepth@##1{\i@=##1\relax\the\Rowdepthtoks@}%
 \colwidthtoks@{\fi}%
 \setbox\zer@\vbox{%
  \unvbox\zer@
  \count@\rowcount@
  \loop
   \unskip\unpenalty
   \setbox\zer@\lastbox
   \ifnum\count@>\maxcolrow@ \advance\count@\m@ne
   \repeat
  \hbox{%
   \unhbox\zer@
   \count@\z@
   \loop
    \unskip
    \setbox\zer@\lastbox
    \edef\next@{\noexpand\or\noexpand\getdim@\the\wd\zer@\the\colwidthtoks@}%
     \global\colwidthtoks@=\expandafter{\next@}%
    \advance\count@\@ne
    \ifnum\count@<\Colcount@
    \repeat}}%
 \edef\next@{\noexpand\ifcase\noexpand\i@\the\colwidthtoks@}%
  \global\colwidthtoks@=\expandafter{\next@}%
 \edef\Colwidth@##1{\i@=##1\relax\the\colwidthtoks@}%
 \colwidthtoks@{}\Rowheighttoks@{}\Rowdepthtoks@{}\widthtoks@{}%
 \Widthtoks@{}\heighttoks@{}\Heighttoks@{}\depthtoks@{}\Depthtoks@{}%
}
\newcount\xoff@
\newcount\yoff@
\newcount\endcount@
\newcount\rcount@
\newdimen\firstx@
\newdimen\firsty@
\newdimen\secondx@
\newdimen\secondy@
\newdimen\tocenter@
\newdimen\charht@
\newdimen\charwd@
\def\outside@{\Err@{This arrow points outside the \string\newCD}}
\newif\ifsvertex@
\newif\iftvertex@
\def\arrow@#1#2{\xoff@=#1\relax\yoff@=#2\relax
 \count@\rowcount@ \advance\count@-\yoff@
 \ifnum\count@<\@ne \outside@ \else \ifnum\count@>\Rowcount@ \outside@ \fi\fi
 \count@\colcount@ \advance\count@\xoff@
 \ifnum\count@<\@ne \outside@ \else \ifnum\count@>\Colcount@ \outside@\fi\fi
 \tcolcount@\colcount@ \advance\tcolcount@\xoff@
 \Width@\rowcount@\colcount@ \tocenter@=-\getdim@ \divide\tocenter@\tw@
 \ifdim\getdim@=\z@
  \firstx@\z@ \firsty@\mathaxis@ \svertex@true
 \else
  \svertex@false
  \ifHshort@
   \Colwidth@\colcount@
    \ifE@ \firstx@=.5\getdim@ \else \firstx@=-.5\getdim@ \fi
  \else
   \ifE@ \firstx@=\getdim@ \else \firstx@=-\getdim@ \fi
   \divide\firstx@\tw@
  \fi
  \ifE@
   \ifH@ \advance\firstx@\thr@@\p@ \else \advance\firstx@-\thr@@\p@ \fi
  \else
   \ifH@ \advance\firstx@-\thr@@\p@ \else \advance\firstx@\thr@@\p@ \fi
  \fi
  \ifN@
   \Height@\rowcount@\colcount@ \firsty@=\getdim@
   \ifV@ \advance\firsty@\thr@@\p@ \fi
  \else
   \ifV@
    \Depth@\rowcount@\colcount@ \firsty@=-\getdim@
    \advance\firsty@-\thr@@\p@
   \else
    \firsty@\z@
   \fi
  \fi
 \fi
 \ifV@
 \else
  \Colwidth@\colcount@
  \ifE@ \secondx@=\getdim@ \else \secondx@=-\getdim@ \fi
  \divide\secondx@\tw@
  \ifE@ \else \getcgap@\colcount@ \advance\secondx@-\getdim@ \fi
  \endcount@=\colcount@ \advance\endcount@\xoff@
  \count@=\colcount@
  \ifE@
   \advance\count@\@ne
   \loop
    \ifnum\count@<\endcount@
    \Colwidth@\count@ \advance\secondx@\getdim@
    \getcgap@\count@ \advance\secondx@\getdim@
    \advance\count@\@ne
    \repeat
  \else
   \advance\count@\m@ne
   \loop
    \ifnum\count@>\endcount@
    \Colwidth@\count@ \advance\secondx@-\getdim@
    \getcgap@\count@ \advance\secondx@-\getdim@
    \advance\count@\m@ne
    \repeat
  \fi
  \Colwidth@\count@ \divide\getdim@\tw@
  \ifHshort@
  \else
   \ifE@ \advance\secondx@\getdim@ \else \advance\secondx@-\getdim@ \fi
  \fi
  \ifE@ \getcgap@\count@ \advance\secondx@\getdim@ \fi
  \rcount@\rowcount@ \advance\rcount@-\yoff@
  \Width@\rcount@\count@ \divide\getdim@\tw@
  \tvertex@false
  \ifH@\ifdim\getdim@=\z@\tvertex@true\Hshort@false\fi\fi
  \ifHshort@
  \else
   \ifE@ \advance\secondx@-\getdim@ \else \advance\secondx@\getdim@ \fi
  \fi
  \iftvertex@
   \advance\secondx@.4\p@
  \else
   \ifE@ \advance\secondx@-\thr@@\p@ \else \advance\secondx@\thr@@\p@ \fi
  \fi
 \fi
 \ifH@
 \else
  \ifN@
   \Rowheight@\rowcount@ \secondy@\getdim@
  \else
   \Rowdepth@\rowcount@ \secondy@-\getdim@
   \getrgap@\rowcount@ \advance\secondy@-\getdim@
  \fi
  \endcount@=\rowcount@ \advance\endcount@-\yoff@
  \count@=\rowcount@
  \ifN@
   \advance\count@\m@ne
   \loop
    \ifnum\count@>\endcount@
    \Rowheight@\count@ \advance\secondy@\getdim@
    \Rowdepth@\count@ \advance\secondy@\getdim@
    \getrgap@\count@ \advance\secondy@\getdim@
    \advance\count@\m@ne
    \repeat
  \else
   \advance\count@\@ne
   \loop
    \ifnum\count@<\endcount@
    \Rowheight@\count@ \advance\secondy@-\getdim@
    \Rowdepth@\count@ \advance\secondy@-\getdim@
    \getrgap@\count@ \advance\secondy@-\getdim@
    \advance\count@\@ne
    \repeat
  \fi
  \tvertex@false
  \ifV@\Width@\count@\colcount@\ifdim\getdim@=\z@\tvertex@true\fi\fi
  \ifN@
   \getrgap@\count@ \advance\secondy@\getdim@
   \Rowdepth@\count@ \advance\secondy@\getdim@
   \iftvertex@
    \advance\secondy@\mathaxis@
   \else
    \Depth@\count@\tcolcount@ \advance\secondy@-\getdim@
    \advance\secondy@-\thr@@\p@
   \fi
  \else
   \Rowheight@\count@ \advance\secondy@-\getdim@
   \iftvertex@
    \advance\secondy@\mathaxis@
   \else
    \Height@\count@\tcolcount@ \advance\secondy@\getdim@
    \advance\secondy@\thr@@\p@
   \fi
  \fi
 \fi
 \ifV@\else\advance\firstx@\sxdimen@\fi
 \ifH@\else\advance\firsty@\sydimen@\fi
 \iftX@
  \advance\secondy@\tXdimen@ii
  \advance\secondx@\tXdimen@i
  \slope@
 \else
  \iftY@
   \advance\secondy@\tYdimen@ii
   \advance\secondx@\tYdimen@i
   \slope@
   \secondy@=\secondx@ \advance\secondy@-\firstx@
   \ifNESW@ \else \multiply\secondy@\m@ne \fi
   \multiply\secondy@\tan@i \divide\secondy@\tan@ii \advance\secondy@\firsty@
  \else
   \ifa@
    \slope@
    \ifNESW@ \global\advance\angcount@\exacount@ \else
      \global\advance\angcount@-\exacount@ \fi
    \ifnum\angcount@>23 \angcount@23 \fi
    \ifnum\angcount@<\@ne \angcount@\@ne \fi
    \slope@a\angcount@
    \ifY@
     \advance\secondy@\Ydimen@
    \else
     \ifX@
      \advance\secondx@\Xdimen@
      \dimen@\secondx@ \advance\dimen@-\firstx@
      \ifNESW@\else\multiply\dimen@\m@ne\fi
      \multiply\dimen@\tan@i \divide\dimen@\tan@ii
      \advance\dimen@\firsty@ \secondy@=\dimen@
     \fi
    \fi
   \else
    \ifH@\else\ifV@\else\slope@\fi\fi
   \fi
  \fi
 \fi
 \ifH@\else\ifV@\else\ifsvertex@\else
  \dimen@=6\p@ \multiply\dimen@\tan@ii
  \count@=\tan@i \advance\count@\tan@ii \divide\dimen@\count@
  \ifE@ \advance\firstx@\dimen@ \else \advance\firstx@-\dimen@ \fi
  \multiply\dimen@\tan@i \divide\dimen@\tan@ii
  \ifN@ \advance\firsty@\dimen@ \else \advance\firsty@-\dimen@ \fi
 \fi\fi\fi
 \ifp@
  \ifH@\else\ifV@\else
   \getcos@\pdimen@ \advance\firsty@\dimen@ \advance\secondy@\dimen@
   \ifNESW@ \advance\firstx@-\dimen@ii \else \advance\firstx@\dimen@ii \fi
  \fi\fi
 \fi
 \ifH@\else\ifV@\else
  \ifnum\tan@i>\tan@ii
   \charht@=10\p@ \charwd@=10\p@
   \multiply\charwd@\tan@ii \divide\charwd@\tan@i
  \else
   \charwd@=10\p@ \charht@=10\p@
   \divide\charht@\tan@ii \multiply\charht@\tan@i
  \fi
  \ifnum\tcount@=\thr@@
   \ifN@ \advance\secondy@-.3\charht@ \else\advance\secondy@.3\charht@ \fi
  \fi
  \ifnum\scount@=\tw@
   \ifE@ \advance\firstx@.3\charht@ \else \advance\firstx@-.3\charht@ \fi
  \fi
  \ifnum\tcount@=12
   \ifN@ \advance\secondy@-\charht@ \else \advance\secondy@\charht@ \fi
  \fi
  \iftY@
  \else
   \ifa@
    \ifX@
    \else
     \secondx@\secondy@ \advance\secondx@-\firsty@
     \ifNESW@\else\multiply\secondx@\m@ne\fi
     \multiply\secondx@\tan@ii \divide\secondx@\tan@i
     \advance\secondx@\firstx@
    \fi
   \fi
  \fi
 \fi\fi
 \ifH@\harrow@\else\ifV@\varrow@\else\arrow@@\fi\fi}
\newdimen\mathaxis@
\mathaxis@90\p@ \divide\mathaxis@36
\def\harrow@b{\ifE@\hskip\tocenter@\hskip\firstx@\fi}
\def\harrow@bb{\ifE@\hskip\xdimen@\else\hskip\Xdimen@\fi}
\def\harrow@e{\ifE@\else\hskip-\firstx@\hskip-\tocenter@\fi}
\def\harrow@ee{\ifE@\hskip-\Xdimen@\else\hskip-\xdimen@\fi}
\def\harrow@{\dimen@\secondx@\advance\dimen@-\firstx@
 \ifE@ \let\next@\rlap \else  \multiply\dimen@\m@ne \let\next@\llap \fi
 \next@{%
  \harrow@b
  \smash{\raise\pdimen@\hbox to\dimen@
   {\harrow@bb\arrow@ii
    \ifnum\arrcount@=\m@ne \else \ifnum\arrcount@=\thr@@ \else
     \ifE@
      \ifnum\scount@=\m@ne
      \else
       \ifcase\scount@\or\or\char118 \or\char117 \or\or\or\char119 \or
       \char120 \or\char121 \or\char122 \or\or\or\arrow@i\char125 \or
       \char117 \hskip\thr@@\p@\char117 \hskip-\thr@@\p@\fi
      \fi
     \else
      \ifnum\tcount@=\m@ne
      \else
       \ifcase\tcount@\char117 \or\or\char117 \or\char118 \or\char119 \or
       \char120\or\or\or\or\or\char121 \or\char122 \or\arrow@i\char125
       \or\char117 \hskip\thr@@\p@\char117 \hskip-\thr@@\p@\fi
      \fi
     \fi
    \fi\fi
    \dimen@\mathaxis@ \advance\dimen@.2\p@
    \dimen@ii\mathaxis@ \advance\dimen@ii-.2\p@
    \ifnum\arrcount@=\m@ne
     \let\leads@\null
    \else
     \ifcase\arrcount@
      \def\leads@{\hrule height\dimen@ depth-\dimen@ii}\or
      \def\leads@{\hrule height\dimen@ depth-\dimen@ii}\or
      \def\leads@{\hbox to10\p@{%
       \leaders\hrule height\dimen@ depth-\dimen@ii\hfil
       \hfil
      \leaders\hrule height\dimen@ depth-\dimen@ii\hskip\z@ plus2fil\relax
       \hfil
       \leaders\hrule height\dimen@ depth-\dimen@ii\hfil}}\or
     \def\leads@{\hbox{\hbox to10\p@{\dimen@\mathaxis@ \advance\dimen@1.2\p@
       \dimen@ii\dimen@ \advance\dimen@ii-.4\p@
       \leaders\hrule height\dimen@ depth-\dimen@ii\hfil}%
       \kern-10\p@
       \hbox to10\p@{\dimen@\mathaxis@ \advance\dimen@-1.2\p@
       \dimen@ii\dimen@ \advance\dimen@ii-.4\p@
       \leaders\hrule height\dimen@ depth-\dimen@ii\hfil}}}\fi
    \fi
    \cleaders\leads@\hfil
    \ifnum\arrcount@=\m@ne\else\ifnum\arrcount@=\thr@@\else
     \arrow@i
     \ifE@
      \ifnum\tcount@=\m@ne
      \else
       \ifcase\tcount@\char119 \or\or\char119 \or\char120 \or\char121 \or
       \char122 \or \or\or\or\or\char123\or\char124 \or
       \char125 \or\char119 \hskip-\thr@@\p@\char119 \hskip\thr@@\p@\fi
      \fi
     \else
      \ifcase\scount@\or\or\char120 \or\char119 \or\or\or\char121 \or\char122
      \or\char123 \or\char124 \or\or\or\char125 \or
      \char119 \hskip-\thr@@\p@\char119 \hskip\thr@@\p@\fi
     \fi
    \fi\fi
    \harrow@ee}}%
  \harrow@e}%
 \iflabel@i
  \dimen@ii\z@ \setbox\zer@\hbox{$\m@th\tsize@@\label@i$}%
  \ifnum\arrcount@=\m@ne
  \else
   \advance\dimen@ii\mathaxis@
   \advance\dimen@ii\dp\zer@ \advance\dimen@ii\tw@\p@
   \ifnum\arrcount@=\thr@@ \advance\dimen@ii\tw@\p@ \fi
  \fi
  \advance\dimen@ii\pdimen@
  \next@{\harrow@b\smash{\raise\dimen@ii\hbox to\dimen@
   {\harrow@bb\hskip\tw@\ldimen@i\hfil\box\zer@\hfil\harrow@ee}}\harrow@e}%
 \fi
 \iflabel@ii
  \ifnum\arrcount@=\m@ne
  \else
   \setbox\zer@\hbox{$\m@th\tsize@\label@ii$}%
   \dimen@ii-\ht\zer@ \advance\dimen@ii-\tw@\p@
   \ifnum\arrcount@=\thr@@ \advance\dimen@ii-\tw@\p@ \fi
   \advance\dimen@ii\mathaxis@ \advance\dimen@ii\pdimen@
   \next@{\harrow@b\smash{\raise\dimen@ii\hbox to\dimen@
    {\harrow@bb\hskip\tw@\ldimen@ii\hfil\box\zer@\hfil\harrow@ee}}\harrow@e}%
  \fi
 \fi}
\let\tsize@\tsize
\def\tsizenewCDlabels{\let\tsize@\tsize}
\def\ssizenewCDlabels{\let\tsize@\ssize}
\def\tsize@@{\ifnum\arrcount@=\m@ne\else\tsize@\fi}
\def\varrow@{\dimen@\secondy@ \advance\dimen@-\firsty@
 \ifN@ \else \multiply\dimen@\m@ne \fi
 \setbox\zer@\vbox to\dimen@
  {\ifN@ \vskip-\Ydimen@ \else \vskip\ydimen@ \fi
   \ifnum\arrcount@=\m@ne\else\ifnum\arrcount@=\thr@@\else
    \hbox{\arrow@iii
     \ifN@
      \ifnum\tcount@=\m@ne
      \else
       \ifcase\tcount@\char117 \or\or\char117 \or\char118 \or\char119 \or
       \char120 \or\or\or\or\or\char121 \or\char122 \or\char123 \or
       \vbox{\hbox{\char117 }\nointerlineskip\vskip\thr@@\p@
       \hbox{\char117 }\vskip-\thr@@\p@}\fi
      \fi
     \else
      \ifcase\scount@\or\or\char118 \or\char117 \or\or\or\char119 \or
      \char120 \or\char121 \or\char122 \or\or\or\char123 \or
      \vbox{\hbox{\char117 }\nointerlineskip\vskip\thr@@\p@
      \hbox{\char117 }\vskip-\thr@@\p@}\fi
     \fi}%
    \nointerlineskip
   \fi\fi
   \ifnum\arrcount@=\m@ne
    \let\leads@\null
   \else
    \ifcase\arrcount@\let\leads@\vrule\or\let\leads@\vrule\or
    \def\leads@{\vbox to10\p@{%
     \hrule height 1.67\p@ depth\z@ width.4\p@
     \vfil
     \hrule height 3.33\p@ depth\z@ width.4\p@
     \vfil
     \hrule height 1.67\p@ depth\z@ width.4\p@}}\or
    \def\leads@{\hbox{\vrule height\p@\hskip\tw@\p@\vrule}}\fi
   \fi
  \cleaders\leads@\vfill\nointerlineskip
   \ifnum\arrcount@=\m@ne\else\ifnum\arrcount@=\thr@@\else
    \hbox{\arrow@iv
     \ifN@
      \ifcase\scount@\or\or\char118 \or\char117 \or\or\or\char119 \or
      \char120 \or\char121 \or\char122 \or\or\or\arrow@iii\char123 \or
      \vbox{\hbox{\char117 }\nointerlineskip\vskip-\thr@@\p@
      \hbox{\char117 }\vskip\thr@@\p@}\fi
     \else
      \ifnum\tcount@=\m@ne
      \else
       \ifcase\tcount@\char117 \or\or\char117 \or\char118 \or\char119 \or
       \char120 \or\or\or\or\or\char121 \or\char122 \or\arrow@iii\char123 \or
       \vbox{\hbox{\char117 }\nointerlineskip\vskip-\thr@@\p@
       \hbox{\char117 }\vskip\thr@@\p@}\fi
      \fi
     \fi}%
   \fi\fi
   \ifN@\vskip\ydimen@\else\vskip-\Ydimen@\fi}%
 \ifN@
  \dimen@ii\firsty@
 \else
  \dimen@ii-\firsty@ \advance\dimen@ii\ht\zer@ \multiply\dimen@ii\m@ne
 \fi
 \rlap{\smash{\hskip\tocenter@ \hskip\pdimen@ \raise\dimen@ii \box\zer@}}%
 \iflabel@i
  \setbox\zer@\vbox to\dimen@{\vfil
   \hbox{$\m@th\tsize@@\label@i$}\vskip\tw@\ldimen@i\vfil}%
  \rlap{\smash{\hskip\tocenter@ \hskip\pdimen@
  \ifnum\arrcount@=\m@ne \let\next@\relax \else \let\next@\llap \fi
  \next@{\raise\dimen@ii\hbox{\ifnum\arrcount@=\m@ne \hskip-.5\wd\zer@ \fi
   \box\zer@ \ifnum\arrcount@=\m@ne \else \hskip\tw@\p@ \fi}}}}%
 \fi
 \iflabel@ii
  \ifnum\arrcount@=\m@ne
  \else
   \setbox\zer@\vbox to\dimen@{\vfil
    \hbox{$\m@th\tsize@\label@ii$}\vskip\tw@\ldimen@ii\vfil}%
   \rlap{\smash{\hskip\tocenter@ \hskip\pdimen@
   \rlap{\raise\dimen@ii\hbox{\ifnum\arrcount@=\thr@@ \hskip4.5\p@ \else
    \hskip2.5\p@ \fi\box\zer@}}}}%
  \fi
 \fi
}
\newdimen\goal@
\newdimen\shifted@
\newcount\Tcount@
\newcount\Scount@
\newbox\shaft@
\newcount\slcount@
\def\getcos@#1{%
 \ifnum\tan@i<\tan@ii
  \dimen@#1%
  \ifnum\slcount@<8 \count@9 \else \ifnum\slcount@<12 \count@8 \else
   \count@7 \fi\fi
  \multiply\dimen@\count@ \divide\dimen@10
  \dimen@ii\dimen@ \multiply\dimen@ii\tan@i \divide\dimen@ii\tan@ii
 \else
  \dimen@ii#1%
  \count@-\slcount@ \advance\count@24
  \ifnum\count@<8 \count@9 \else \ifnum\count@<12 \count@8
   \else\count@7 \fi\fi
  \multiply\dimen@ii\count@ \divide\dimen@ii10
  \dimen@\dimen@ii \multiply\dimen@\tan@ii \divide\dimen@\tan@i
 \fi}
\newdimen\adjust@
\def\Nnext@{\ifN@\let\next@\raise\else\let\next@\lower\fi}
\def\arrow@@{\slcount@\angcount@
 \ifNESW@
  \ifnum\angcount@<10
   \let\arrowfont@=\arrow@i \advance\angcount@\m@ne \multiply\angcount@13
  \else
   \ifnum\angcount@<19
    \let\arrowfont@=\arrow@ii \advance\angcount@-10 \multiply\angcount@13
   \else
    \let\arrowfont@=\arrow@iii \advance\angcount@-19 \multiply\angcount@13
  \fi\fi
  \Tcount@\angcount@
 \else
  \ifnum\angcount@<5
   \let\arrowfont@=\arrow@iii \advance\angcount@\m@ne \multiply\angcount@13
   \advance\angcount@65
  \else
   \ifnum\angcount@<14
    \let\arrowfont@=\arrow@iv \advance\angcount@-5 \multiply\angcount@13
   \else
    \ifnum\angcount@<23
     \let\arrowfont@=\arrow@v \advance\angcount@-14 \multiply\angcount@13
    \else
     \let\arrowfont@=\arrow@i \angcount@=117
  \fi\fi\fi
  \ifnum\angcount@=117 \Tcount@=115 \else\Tcount@\angcount@ \fi
 \fi
 \Scount@\Tcount@
 \ifE@
  \ifnum\tcount@=\z@ \advance\Tcount@\tw@ \else\ifnum\tcount@=13
   \advance\Tcount@\tw@ \else \advance\Tcount@\tcount@ \fi\fi
  \ifnum\scount@=\z@ \else \ifnum\scount@=13 \advance\Scount@\thr@@ \else
   \advance\Scount@\scount@ \fi\fi
 \else
  \ifcase\tcount@\advance\Tcount@\thr@@\or\or\advance\Tcount@\thr@@\or
  \advance\Tcount@\tw@\or\advance\Tcount@6 \or\advance\Tcount@7
  \or\or\or\or\or \advance\Tcount@8 \or\advance\Tcount@9 \or
  \advance\Tcount@12 \or\advance\Tcount@\thr@@\fi
  \ifcase\scount@\or\or\advance\Scount@\thr@@\or\advance\Scount@\tw@\or
  \or\or\advance\Scount@4 \or\advance\Scount@5 \or\advance\Scount@10
  \or\advance\Scount@11 \or\or\or\advance\Scount@12 \or\advance
  \Scount@\tw@\fi
 \fi
 \ifcase\arrcount@\or\or\advance\angcount@\@ne\else\fi
 \ifN@ \shifted@=\firsty@ \else\shifted@=-\firsty@ \fi
 \ifE@ \else\advance\shifted@\charht@ \fi
 \goal@=\secondy@ \advance\goal@-\firsty@
 \ifN@\else\multiply\goal@\m@ne\fi
 \setbox\shaft@\hbox{\arrowfont@\char\angcount@}%
 \ifnum\arrcount@=\thr@@
  \getcos@{1.5\p@}%
  \setbox\shaft@\hbox to\wd\shaft@{\arrowfont@
   \rlap{\hskip\dimen@ii
    \smash{\ifNESW@\let\next@\lower\else\let\next@\raise\fi
     \next@\dimen@\hbox{\arrowfont@\char\angcount@}}}%
   \rlap{\hskip-\dimen@ii
    \smash{\ifNESW@\let\next@\raise\else\let\next@\lower\fi
      \next@\dimen@\hbox{\arrowfont@\char\angcount@}}}\hfil}%
 \fi
 \rlap{\smash{\hskip\tocenter@\hskip\firstx@
  \ifnum\arrcount@=\m@ne
  \else
   \ifnum\arrcount@=\thr@@
   \else
    \ifnum\scount@=\m@ne
    \else
     \ifnum\scount@=\z@
     \else
      \setbox\zer@\hbox{\ifnum\angcount@=117 \arrow@v\else\arrowfont@\fi
       \char\Scount@}%
      \ifNESW@
       \ifnum\scount@=\tw@
        \dimen@=\shifted@ \advance\dimen@-\charht@
        \ifN@\hskip-\wd\zer@\fi
        \Nnext@
        \next@\dimen@\copy\zer@
        \ifN@\else\hskip-\wd\zer@\fi
       \else
        \Nnext@
        \ifN@\else\hskip-\wd\zer@\fi
        \next@\shifted@\copy\zer@
        \ifN@\hskip-\wd\zer@\fi
       \fi
       \ifnum\scount@=12
        \advance\shifted@\charht@ \advance\goal@-\charht@
        \ifN@ \hskip\wd\zer@ \else \hskip-\wd\zer@ \fi
       \fi
       \ifnum\scount@=13
        \getcos@{\thr@@\p@}%
        \ifN@ \hskip\dimen@ \else \hskip-\wd\zer@ \hskip-\dimen@ \fi
        \adjust@\shifted@ \advance\adjust@\dimen@ii
        \Nnext@
        \next@\adjust@\copy\zer@
        \ifN@ \hskip-\dimen@ \hskip-\wd\zer@ \else \hskip\dimen@ \fi
       \fi
      \else
       \ifN@\hskip-\wd\zer@\fi
       \ifnum\scount@=\tw@
        \ifN@ \hskip\wd\zer@ \else \hskip-\wd\zer@ \fi
        \dimen@=\shifted@ \advance\dimen@-\charht@
        \Nnext@
        \next@\dimen@\copy\zer@
        \ifN@\hskip-\wd\zer@\fi
       \else
        \Nnext@
        \next@\shifted@\copy\zer@
        \ifN@\else\hskip-\wd\zer@\fi
       \fi
       \ifnum\scount@=12
        \advance\shifted@\charht@ \advance\goal@-\charht@
        \ifN@ \hskip-\wd\zer@ \else \hskip\wd\zer@ \fi
       \fi
       \ifnum\scount@=13
        \getcos@{\thr@@\p@}%
        \ifN@ \hskip-\wd\zer@ \hskip-\dimen@ \else \hskip\dimen@ \fi
        \adjust@\shifted@ \advance\adjust@\dimen@ii
        \Nnext@
        \next@\adjust@\copy\zer@
        \ifN@ \hskip\dimen@ \else \hskip-\dimen@ \hskip-\wd\zer@ \fi
       \fi	
      \fi
  \fi\fi\fi\fi
  \ifnum\arrcount@=\m@ne
  \else
   \loop
    \ifdim\goal@>\charht@
    \ifE@\else\hskip-\charwd@\fi
    \Nnext@
    \next@\shifted@\copy\shaft@
    \ifE@\else\hskip-\charwd@\fi
    \advance\shifted@\charht@ \advance\goal@ -\charht@
    \repeat
   \ifdim\goal@>\z@
    \dimen@=\charht@ \advance\dimen@-\goal@
    \divide\dimen@\tan@i \multiply\dimen@\tan@ii
    \ifE@ \hskip-\dimen@ \else \hskip-\charwd@ \hskip\dimen@ \fi
    \adjust@=\shifted@ \advance\adjust@-\charht@ \advance\adjust@\goal@
    \Nnext@
    \next@\adjust@\copy\shaft@
    \ifE@ \else \hskip-\charwd@ \fi
   \else
    \adjust@=\shifted@ \advance\adjust@-\charht@
   \fi
  \fi
  \ifnum\arrcount@=\m@ne
  \else
   \ifnum\arrcount@=\thr@@
   \else
    \ifnum\tcount@=\m@ne
    \else
     \setbox\zer@
      \hbox{\ifnum\angcount@=117 \arrow@v\else\arrowfont@\fi\char\Tcount@}%
     \ifnum\tcount@=\thr@@
      \advance\adjust@\charht@
      \ifE@\else\ifN@\hskip-\charwd@\else\hskip-\wd\zer@\fi\fi
     \else
      \ifnum\tcount@=12
       \advance\adjust@\charht@
       \ifE@\else\ifN@\hskip-\charwd@\else\hskip-\wd\zer@\fi\fi
      \else
       \ifE@\hskip-\wd\zer@\fi
     \fi\fi
     \Nnext@
     \next@\adjust@\copy\zer@
     \ifnum\tcount@=13
      \hskip-\wd\zer@
      \getcos@{\thr@@\p@}%
      \ifE@\hskip-\dimen@ \else\hskip\dimen@ \fi
      \advance\adjust@-\dimen@ii
      \Nnext@
      \next@\adjust@\box\zer@
     \fi
  \fi\fi\fi}}%
 \iflabel@i
  \rlap{\hskip\tocenter@
  \dimen@\firstx@ \advance\dimen@\secondx@ \divide\dimen@\tw@
  \advance\dimen@\ldimen@i
  \dimen@ii\firsty@ \advance\dimen@ii\secondy@ \divide\dimen@ii\tw@
  \multiply\ldimen@i\tan@i \divide\ldimen@i\tan@ii
  \ifNESW@ \advance\dimen@ii\ldimen@i \else \advance\dimen@ii-\ldimen@i \fi
  \setbox\zer@\hbox{\ifNESW@\else\ifnum\arrcount@=\thr@@\hskip4\p@\else
   \hskip\tw@\p@\fi\fi
   $\m@th\tsize@@\label@i$\ifNESW@\ifnum\arrcount@=\thr@@\hskip4\p@\else
   \hskip\tw@\p@\fi\fi}%
  \ifnum\arrcount@=\m@ne
   \ifNESW@ \advance\dimen@.5\wd\zer@ \advance\dimen@\p@ \else
    \advance\dimen@-.5\wd\zer@ \advance\dimen@-\p@ \fi
   \advance\dimen@ii-.5\ht\zer@
  \else
   \advance\dimen@ii\dp\zer@
   \ifnum\slcount@<6 \advance\dimen@ii\tw@\p@ \fi
  \fi
  \hskip\dimen@
  \ifNESW@ \let\next@\llap \else\let\next@\rlap \fi
  \next@{\smash{\raise\dimen@ii\box\zer@}}}%
 \fi
 \iflabel@ii
  \ifnum\arrcount@=\m@ne
  \else
   \rlap{\hskip\tocenter@
   \dimen@\firstx@ \advance\dimen@\secondx@ \divide\dimen@\tw@
   \ifNESW@ \advance\dimen@\ldimen@ii \else \advance\dimen@-\ldimen@ii \fi
   \dimen@ii\firsty@ \advance\dimen@ii\secondy@ \divide\dimen@ii\tw@
   \multiply\ldimen@ii\tan@i \divide\ldimen@ii\tan@ii
   \advance\dimen@ii\ldimen@ii
   \setbox\zer@\hbox{\ifNESW@\ifnum\arrcount@=\thr@@\hskip4\p@\else
    \hskip\tw@\p@\fi\fi
    $\m@th\tsize@\label@ii$\ifNESW@\else\ifnum\arrcount@=\thr@@\hskip4\p@
    \else\hskip\tw@\p@\fi\fi}%
   \advance\dimen@ii-\ht\zer@
   \ifnum\slcount@<9 \advance\dimen@ii-\thr@@\p@ \fi
   \ifNESW@ \let\next@\rlap \else \let\next@\llap \fi
   \hskip\dimen@\next@{\smash{\raise\dimen@ii\box\zer@}}}%
  \fi
 \fi
}
\def\outnewCD@#1{\def#1{\Err@{\string#1 must not be used within \string\newCD}}}
\newskip\prenewCDskip@
\newskip\postnewCDskip@
\prenewCDskip@\z@
\postnewCDskip@\z@
\def\prenewCDspace#1{\RIfMIfI@
 \onlydmatherr@\prenewCDspace\else\advance\prenewCDskip@#1\relax\fi\else
 \onlydmatherr@\prenewCDspace\fi}
\def\postnewCDspace#1{\RIfMIfI@
 \onlydmatherr@\postnewCDspace\else\advance\postnewCDskip@#1\relax\fi\else
 \onlydmatherr@\postnewCDspace\fi}
\def\predisplayspace#1{\RIfMIfI@
 \onlydmatherr@\predisplayspace\else
 \advance\abovedisplayskip#1\relax
 \advance\abovedisplayshortskip#1\relax\fi
 \else\onlydmatherr@\prenewCDspace\fi}
\def\postdisplayspace#1{\RIfMIfI@
 \onlydmatherr@\postdisplayspace\else
 \advance\belowdisplayskip#1\relax
 \advance\belowdisplayshortskip#1\relax\fi
 \else\onlydmatherr@\postdisplayspace\fi}
\def\PrenewCDSpace#1{\global\prenewCDskip@#1\relax}
\def\PostnewCDSpace#1{\global\postnewCDskip@#1\relax}
\def\newCD#1\endnewCD{%
 \outnewCD@\cgaps\outnewCD@\rgaps\outnewCD@\Cgaps\outnewCD@\Rgaps
 \prenewCD@#1\endnewCD
 \advance\abovedisplayskip\prenewCDskip@
 \advance\abovedisplayshortskip\prenewCDskip@
 \advance\belowdisplayskip\postnewCDskip@
 \advance\belowdisplayshortskip\postnewCDskip@
 \vcenter{\vskip\prenewCDskip@ \Let@ \colcount@\@ne \rowcount@\z@
  \everycr{%
   \noalign{%
    \ifnum\rowcount@=\Rowcount@
    \else
     \global\nointerlineskip
     \getrgap@\rowcount@ \vskip\getdim@
     \global\advance\rowcount@\@ne \global\colcount@\@ne
    \fi}}%
  \tabskip\z@
  \halign{&\global\xoff@\z@ \global\yoff@\z@
   \getcgap@\colcount@ \hskip\getdim@
   \hfil\vrule height10\p@ width\z@ depth\z@
   $\m@th\displaystyle{##}$\hfil
   \global\advance\colcount@\@ne\cr
   #1\crcr}\vskip\postnewCDskip@}%
 \prenewCDskip@\z@\postnewCDskip@\z@
 \def\getcgap@##1{\ifcase##1\or\getdim@\z@\else\getdim@\standardcgap\fi}%
 \def\getrgap@##1{\ifcase##1\getdim@\z@\else\getdim@\standardrgap\fi}%
 \let\Width@\relax\let\Height@\relax\let\Depth@\relax\let\Rowheight@\relax
 \let\Rowdepth@\relax\let\Colwdith@\relax
}
\catcode`\@=\active